%% file: the_dirichlet_random_walk.tex
\documentclass[10pt,a4paper]{amsart}
\input{preambule_spectral_gap}

\title{The Dirichlet random walk.}

\begin{document}

\author{Adrien Boulanger}

\author{Olivier Glorieux}

\thanks{This project received funding from the European Research Council (ERC) under the European Union's Horizon 2020 research and innovation programme (ERC 647133 IChaos and ERC 715982 DiGGeS)}

\maketitle

\begin{abstract}
In this article we define and study a stochastic process on Galoisian covers of compact manifolds. The successive positions of the process are defined recursively by picking a point uniformly in the Dirichlet domain of the previous one. We prove a theorem \textit{à la Kesten} for such a process: the escape rate of the random walk is positive if and only if the cover is non amenable. We also investigate more in details the case where the deck group is Gromov hyperbolic, showing the almost sure convergence to the boundary of the trajectory as well as a central limit theorem for the escape rate. 
\end{abstract}

\section{Introduction} 

Let $M$ be a complete connected Riemannian manifold of dimension $d$ and $ \pi : M \to M_0 := \quotient{M}{\Gamma}$ be a Galoisian Riemannian covering of deck group $\Gamma$ with $M_0$ compact. \\ 

This article aims at studying the large time behavior of a stochastic process on $M$ which is constructed as follows. Let $o \in M$ be any given point and $p_0$ its image in $M_0$. Pick independently and uniformly with respect to the Riemannian measure countably many points $(X_i)_{i \in \NN^*}$ in the compact manifold $M_0$. For any $i \ge 0$, we denote by $g_i$ the (almost surely unique) minimizing geodesic whose endpoints are $X_i$ and $X_{i+1}$. Concatenating the $g_i$s altogether provides us with a random piecewise geodesic path of $M_0$. Such a path being in particular continuous,  it can be lifted on $M$ as a path starting at $o$. We denote by $Z_n \in M $  the endpoint of the piecewise geodesic which corresponds to the concatenation of $g_1,..., g_n$. By an abuse of notations, we denote $Z_0 = o$. The resulting stochastic process is called the Dirichlet random walk. Note in particular that the above construction makes sense with the universal cover of $M_0$.

\begin{figure}[h!]
\begin{center}
	\def\svgwidth{0.7 \columnwidth}
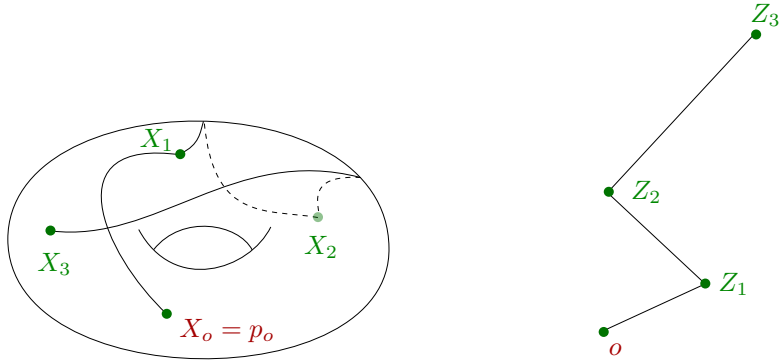
	\end{center}
\caption{On the left, the compact quotient manifold $M_0$. The random points $X_i$s, in green, are related by the geodesics $g_i$, in black. The concatenated curve is then lift to the cover $M$, represented on the right.}
	\label{fig.illusdirichletrandomwalk}
\end{figure}

The study of this stochastic process is motivated in part by its geometric flavour: the behaviour of the process depends on the geometry of both the compact manifold and the covering group. For example in the case where $M_0$ is a hyperbolic surface and $M = \HH^2$ its universal cover, the process depends on the action of $\pi_1(M_0)$ on $\HH^2$ (on the contrary of the Brownian motion for instance). We will see that this process can be completely described with the so called Dirichlet domains that we defined in Section \ref{sec -overview}: knowing that the process is at $x$ at some time, the law of the next point will be given by the characteristic function of the Dirichlet domains associated to $x$. These domains depends on $x$, which is why one cannot expect in general to realise a Dirichlet random walk as the pushforward of a random walk on the isometry group, that we call standard random walks. \\

However exceptions exist, like the Dirichlet random walks associated to flat tori and their covers. Indeed we show, in the appendix of this article, that in this case the Dirichlet random walk behaves like a summation of I.I.Ds random variables. Also, in the appendix, we will investigate in general when a Dirichlet random walk can be described as a standard random walk, showing that the flat torus case being the essentially only non compact case. \\

From far away, the geometry of the cover $M$ looks like the one of the deck group endowed with any word metric (the deck groups is finitely generated by Milnor-Svac lemma). In particular these two metric spaces are roughly isometric as defined in \cite{artkanai} (in particular quasi-isometric). One might then expect the  Dirichlet random walk to roughly behave like a symmetric finitely generated one on a Cayley graph of the deck group $\Gamma$. \\

In this article, we study the long time behaviour of such a stochastic process with a special focus on the Gromov hyperbolic case. We refer to \cite{GdlH} for generalities on Gromov hyperbolicity. Recall that $M$ is Gromov hyperbolic if and only if the deck group $\G$ is also Gromov hyperbolic (they are both length spaces and quasi-isometric to one another). In particular, the Gromov boundary of $M$ is homeomorphic to the one of $\Gamma$. The main result of the paper is the following statement.
 
\begin{theorem}\label{th-convergence dans le bord intro}
Suppose the deck group $\G$ is non amenable then there exists $\ell>0$ such that the following convergence holds almost surely:
$$\lim_{n\tv \infty} \frac{d(o,Z_n)}{n}=\ell.$$
If moreover $\G$ is Gromov hyperbolic then $Z_n$ converges almost surely to a point in the Gromov boundary of $M$.
\end{theorem}

Then using results of Mathieu-Sisto \cite{artmathieusisto} we obtain a Central Limit Theorem: 
\begin{theorem}
\label{th-main-asymptotic-behaviour-of-Zn}
Suppose the deck group $\G$ is non-amenable and Gromov hyperbolic then the sequence $(d(o, Z_n))_{n \in \NN}$ satisfies a central limit theorem. Namely, the sequence of random variables
	 $$ \left( \frac{d(o, Z_n) - l n}{\sqrt{n}} \right)_{n \in \NN}  $$
converges in law to a Gaussian random variable.
\end{theorem}

The above results are more or less classical for random walks on hyperbolic groups or for the Brownian motion on the universal cover of compact negatively curved manifolds. For the random walk aspect, the central limit theorem was proven in \cite{artscandcentrallimit} under a finite exponential moment  and in \cite{artbqcentrallimit} under a finite second moment. For the Brownian motion on the universal cover of a compact negatively curved manifold the central limit theorem was proven in \cite{ledrappierCLT}. \\

Following the classical line of work, in order to study the behaviour of $Z_n$, we look to the spectral properties of the Markov  operator associated to this random process, that we name Dirichlet operator, see section \ref{sec -overview}. We obtain a Theorem \emph{à la Kesten} namely,
\begin{theorem}[spectral gap]
	\label{theo-Kesten for Dirichlet}
		\label{theo-spectral gap}
	Let $M_0$ be a close  Riemannian manifold and $M \to M_0$ be a Galoisian Riemannian covering of deck group $\Gamma$. 
The group $\G$ is not amenable if and only if the Dirichlet operator has a spectral gap. 
\end{theorem}

{As already emphasised,} the Dirichlet random walk does not come from any pushforward of a random walk on a group and as such does not fit in the range of Kesten's criterion \cite{artkestenamenable} for non amenability. Note also that the Dirichlet random walk is not the time $1$ of a diffusion (the transition kernel being not even continuous) and as such does not fit in the range of application of the theory developed for Brownian motion. \\

Many results which are well known in those cases could be interesting to investigate in the case of the Dirichlet random walk, like the large deviations theory or whether or not the Martin boundary identifies with the Gromov boundary in the case where the deck group is hyperbolic. The study of the corresponding hitting measure on the Gromov boundary  also seems interesting to the authors. For example being given two hyperbolic metrics on a given topological surfaces  is there something that can be said on whether or not the corresponding hitting measures are singular to one another ? More precisely, if $S$ is a closed surface of genus $\ge 2$ then for any hyperbolic metric (that we think of as a representation $\rho : \pi_1(S) \to \Isom_+(\HH^2)$) one gets both a hitting measure $\nu_\rho$ and an identification $\Phi_{\rho} : \partial \pi_1(S) \to \partial \HH^2$ (where $\partial X$ denotes the Gromov boundary of the hyperbolic space $X$). Does the class of the measure $\Phi_{\rho}^*(\nu_{\rho})$ depend on $\rho$ ? \\

\textbf{Acknowledgements.} The first author wants to thank Antoine Julia and Nicola Tholozan. The first one for his useful explanations on the notion of perimeter and how it relates with the co-area formula and the second one for helping to prove Proposition \ref{prop.dirichletegalstandard}. Both the authors want to thank Gilles Courtois for useful comments on this work. The second author want to thank many supportive persons, Itai Benjamini, Pierre-Louis Blayac, Peter Haissinski and  François Maucourant among others.

\section{Overview of the article.} \label{sec -overview}

This section is devoted to introducing the different objects of this article and sketching the proofs of the above theorems. \\

Let $M$ be a Riemannian manifold and $ \pi : M \to M_0 := \quotient{M}{\Gamma}$ a Galoisian Riemannian covering of deck group $\Gamma$ with $(M_0,g)$ close. Since $M$ and $M_0$ are locally isometric we will denote, by an abuse of notations, both the metrics by $g$ and both the Riemannian measures by $\mu_g$. \\

Let $(\Omega, \PP)$ be a probability space and $(X_i)_{i \in \NN}$ be countably many I.I.D. random variables from $\Omega$ to $M_0$ which follow the law of $\mu_g$ normalized to have total mass $1$. We refer to such a sequence as the \textbf{increments} of the random walk. \\

\begin{figure}[h!]
\begin{center}
	\def\svgwidth{0.7 \columnwidth}
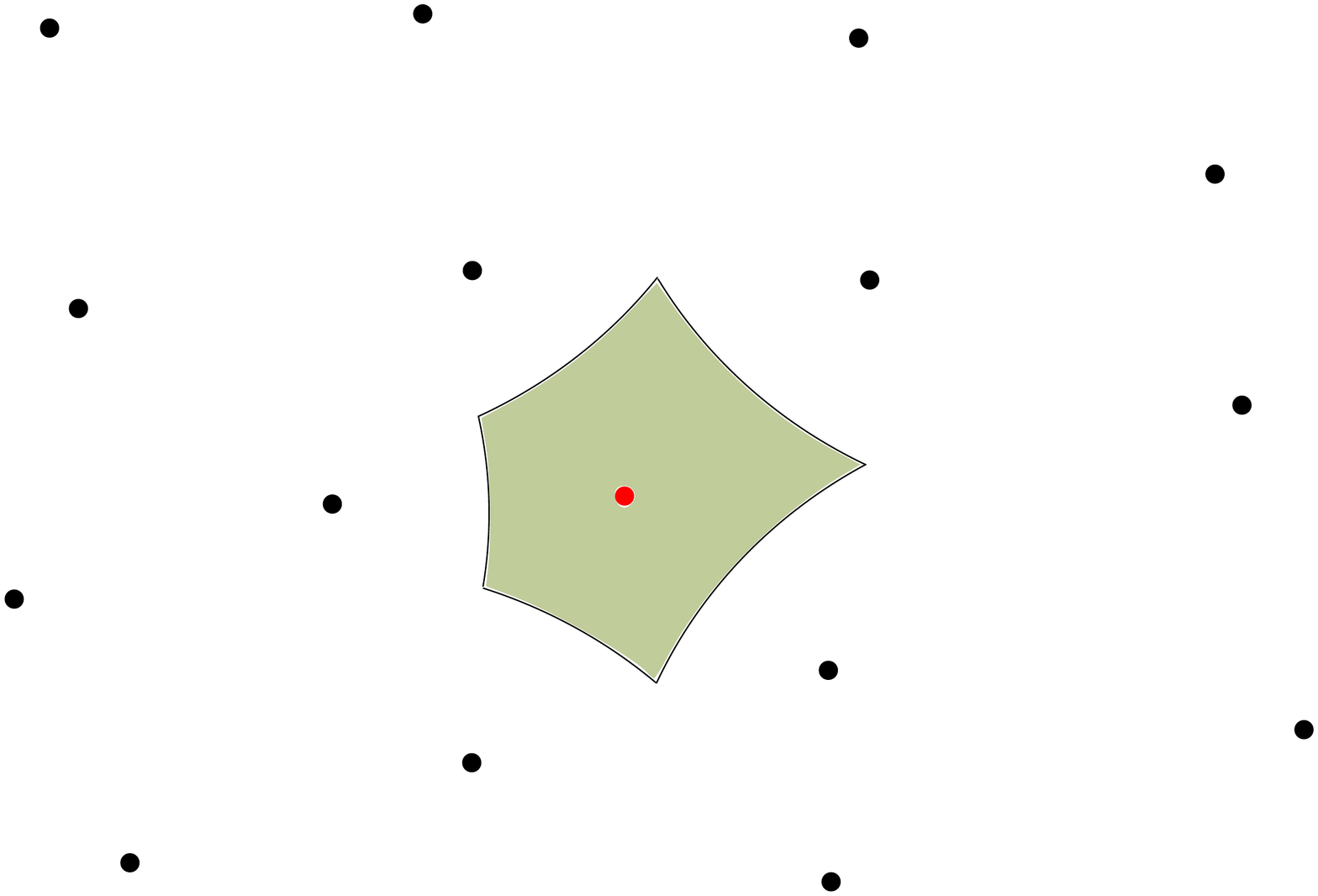
	\end{center}
\caption{The Dirichlet domain of $x$, in green here, can be thought of a polygon as it defined as the intersection of all half spaces containing $x$ delimited by the median planes between $x$ and $\gamma \cdot x$. As The Dirichlet domain is in particular a fundamental domain, the elements of $\Gamma$ are pairing its sides, the median spaces, together. The topological space resulting from these pairings is nothing but the compact quotient manifold $M_0$.}
	\label{fig.dirichletdomain}
\end{figure}

A central notion in this work is the one of  Dirichlet domains, illustrated in Figure \ref{fig.dirichletdomain}.

\begin{definition}
	\label{def-dirichlet domain}
If $X$ is a metric space and $\Gamma$ is a subgroup of $\Isom(X)$ acting properly and discontinuously, we define the \textbf{Dirichlet domain} centred at $x\in X$, denoted $D_x$, as the subset of $M$ defined by
$$D_x :=\{ z \in M \ | \ d( z, x) < d(z , \gamma \cdot x) \ , \  \forall \gamma \in \Gamma \setminus \Id \}. $$
\end{definition}

The cut locus of a Riemannian manifold being always of zero measure (see for example \cite{artcutlocus}) $D_x$ is canonically measurably equivalent to $M_0$ for any $x \in X$.  The main property of Dirichlet domain is that for all $y\in D_x$, 
$d_M(x,y) = d_{M_0}(\pi(x) ,\pi(y))$. Therefore the minimizing geodesics on $D_x$ between $y$ and $x$ projects to the minimizing geodesic on $M_0$. Note also that Dirichlet domains are symmetric, that is $y \in D_x$ if and only if $x \in D_y$. \\

Let us now see how to describe the random process described in the introduction through lifts in Dirichlet domains. Choose a starting point $o = Z_0 \in M$ and pick recursively a random point $Z_n$ in $D_{Z_{n-1}} $ with respect to the normalised Riemannian measure. More formally, one defines reccursively $Z_n$ as follows:
$$Z_n : \Omega \tv M$$
$$\omega \mapsto D_{Z_{n-1}(\omega)} (X_n(\omega))), $$
where we denoted by $ D_x(y)$ the unique lift of a point $y \in M_0$ in the Dirichlet domain $D_x$, $x \in M$. \\

We refer to the random variables $Z_n$ as the \textbf{positions} of the random walk. The resulting stochastic process is called the \textbf{Dirichlet random walk}. Note that the position at time $n$ is completely determined by the $n$ first increments. Note moreover, by independence of the $(X_i)$, that the Dirichlet random walk is a Markov chain on $M$ with transition kernel given by  
	$$ p(x,y) := \frac{1}{\mu_g(M_0)} \mathds{1}_{D_x}(y) \ , $$ 
with respect to the Riemannian measure $\mu_g$. Note that $p(x,y) = p(y,x)$ since Dirichlet domains are symmetric. \\ 

 The study of the random process $Z_n$ relies on the behaviour of the Markov operator associated to the transition kernel $p$. We call this operator the Dirichlet operator. We denote by $L^2(M)$ the set of $\mu_g$ square integrable functions. For any $f_1,f_2 \in L^2(M)$ we denote by 
	$$ \left< f_1,f_2 \right> := \inte{M} f_1 f_2 \ d \mu_g $$
the standard scalar product on $L^2(M)$ and by $|| \cdot ||_{2}$ the associate norm. 

\begin{definition}
	For any $\epsilon > 0$ we call \textbf{Dirichlet operator}, which we denote by $\mathcal{O}$, the operator acting on $L^2(M)$ defined as 
	$$ \mathcal{O}(f)(x) := \inte{M} \  p(x,y)f(y) \ d \mu_g(y) \ .$$
\end{definition}

As any Markov operator, the Dirichlet operator has spectral radius less or equal than one and we denote by $|| \mathcal{O} ||$ its operator norm. Moreover, since $p(x,y)$ is symmetric, the operator $\mathcal{O}$ is self-adjoint with respect to $\left<\cdot , \cdot \right>$. \\

The cornerstone of this article is Theorem \ref{theo-Kesten for Dirichlet} that asserts that the group  $\G$ is not amenable if and only if $|| \mathcal{O}|| < 1$. \\

For the necessary part, the proof is straightforward as in Kesten's proof: we exhibit a  sequence $(f_n)$ of normalised $L^2$ functions for which the norm of $(\mathcal{O}(f_n))$ tends to $1$ as $n \to \infty$. The construction relies on F\o lner criterion of amenability. The proof of the sufficient part, the technical heart of this article, follows the steps of the proof of Cheeger's inequality. Ultimately, it reduces to the proof of a so-called  \textbf{non-local isoperimetric inequality}: 
\begin{proposition}
\label{prop-isoperimetric inequality} 
Under the assumptions of Theorem \ref{theo-spectral gap}, there exists $\beta >0$ such that for all relatively compact open sets $U\subset M$:
$$ \inte{U} \inte{U^c} p^{*2}(x,y) \ d \mu_g(y) d \mu_g(x) \geq \beta \ \mu_g (U) \ , $$ 
where $p^{*2}(x,y)$ denotes the kernel of the operator $\mathcal{O}^2 := \mathcal{O} \circ \mathcal{O}$.
\end{proposition}
We refer to the survey \cite{surveynonlocalperimeter} for more details on non local isoperimetric inequalities and to \cite{artnonlocalperimeter} for a similar statement in a different setting. 

\begin{remark}
The proof of the above proposition does not require to work either with a covering of a compact manifold nor with a Markov operator. One can adapt the proof to the setting of a complete manifold with Ricci curvature and injectivity radius bounded from below and for symmetric positive kernels with the property that there is $\epsilon, C > 0$ such that for any $x \in M$ 
	$$ \inte{B(x, \epsilon)} p(x,y) \ d \mu_g(y) \ge C \ . $$
	In particular, Proposition \ref{prop-isoperimetric inequality} holds for kernels of the form $p(x,y)=  \mathds{1}_{B(x,\epsilon)} +\eta(x,y)$, where $\eta$ is \emph{any} positive symmetric kernel. 
However, to get the probabilist interpretation one would like to require moreover that for any $y \in M$
	$$\inte{M} p(x,y) \ d \mu_g(x) =1 \ ,$$
in other words, that the kernel is Markov. There is actually not so many natural examples. Note for example that walking by drawing uniformly a point on a ball of radius $1$ around where you are does not give rise to a symmetric kernel $p(x,y)$ (unless all the balls of radius $1$ have same volume). 
\end{remark}

Besides the characterization of Theorem \ref{theo-spectral gap},  we will use the spectral gap for non-amenable covers to prove Theorem \ref{th-convergence dans le bord intro}. Indeed, a classical argument, that we will recall in Subsection \ref{subsec.linearprogresswithexp}, shows that the spectral gap  implies the \textbf{linear progress with exponential tail property} of the walk. Namely, the existence of $\epsilon, c > 0$ such that for any $ n,m \in \NN$ we have
\begin{equation}
	\label{eqdef-linearprogress}
	 \PP( d(Z_n,Z_m) \leq \epsilon |n-m|) \le c^{-1} e^{-c |n-m|} \ .
\end{equation}
It follows easily that the Dirichlet random walk is transcient and converges almost surely in the Gromov boundary of the the cover $M$. \\
 
The proof of Theorem \ref{th-main-asymptotic-behaviour-of-Zn} relies on the material developed  in \cite{artmathieusisto}: it reduces to proving that a certain random variable has a finite second moment. In order to introduce this random variable, we recall the definition of the Gromov product. If $X$ is a metric space and $x,y,o \in X$ we define the Gromov product of $x,y$ seen from $o$ as 
	$$ \left<x,y \right>_o := \frac{1}{2} \left( d(x,o) + d(y,o) - d(x,y) \right) \ .$$

Note that for a geodesic $\delta$-hyperbolic space $X$, the quantity $\left<x,y \right>_0$ can be interpreted, up to an additive constant that only depends on $\delta$ as the distance between $o$ and any geodesics from $x$ to $y$. The main difficulty that one encounters in order to use the content of \cite{artmathieusisto} is that one has to show that the walk satisfies to a so-called second moment deviation inequality. In fact, we will prove \textbf{an exponential moment deviation inequality}. Namely, that there is $\epsilon > 0$ such that for any $ n \in \NN$, any $0 \le k \le n$ and any $t > 0$ we have
$$\bP(\left<Z_n,Z_0 \right>_{Z_k}  \geq  t ) = \epsilon^{-1} e^{-\epsilon t} \ . $$

We prove the deviation inequality in Proposition in the appendix of this article, as the proof follows the line of the one proposed in \cite{artmathieusisto}. We decided to give here the full proof to emphasis that no assumption of independence is required whatsoever (on the contrary of how it is stated in \cite{artmathieusisto} or in \cite{artpierremoilargedev}).

\begin{proposition}
\label{prop.morseprobabilistic}
	Let $(Z_n)_{n \in \NN}$ be a sequence of random variables valued in a geodesic Gromov hyperbolic space such that:
		\begin{itemize}
			\item there is $R > 0$ such that for any $n \in \NN$, $d(Z_n, Z_{n+1}) \le R$
			\item the sequence $(Z_n)_{n \in \NN}$ satisfies the linear progress with exponential tail property.
		\end{itemize}
	Then it satisfies an exponential deviation inequality.
\end{proposition}

Note in particular that the above proposition applies in the case of the successive positions of a \textbf{deterministic} quasi-geodesic: the conclusion is then exactly the content of Morse's Lemma (all probabilities are $0$ or $1$ in the deterministic case). \\

Being more careful, following the argumentation of \cite[Proposition 8.2]{artpierremoilargedev}, one can actually remove the assumption that the walk has bounded jumps for the one of uniform finite exponential moment for the law of the jumps for the same conclusions. One can also assume a finite polynomial moment on the law of the jump and this will lead to finite polynomial moment of the deviation inequality. This will have to follow the original approach of \cite[Section 11]{artmathieusisto}.

\section{Spectral gap Theorem}

\label{sec-proof spectral gap}

\subsection{The necessary direction}
\label{subsec-if and only if}
	This section is dedicated to proving that if $\Gamma$ is amenable then the Dirichlet operator $\mathcal{O}$ has no spectral gap. The proof follows a classical strategy: using the F\o lner criterion, we construct a sequence of $L^2$ functions, for which we can bound the operator norm. \\
	
\textbf{Proof of ($\Gamma$  amenable $\Longrightarrow$ no spectral gap).} Since we assumed that $\Gamma$ acts co-compactly on $M$, Milnor-Svarc lemma implies that $\Gamma$ is finitely generated.  We fix a generating system and we identify  $\Gamma$ to the metric space given by the corresponding Cayley graph. Endowed with such a distance, it is well known that $\Gamma$ is quasi-isometric to $M$. \\

For any $\Omega \subset \Gamma$, we denote by $\partial \Omega$ the set of all edges such that one of its endpoints lies in $\Omega$ and the other one lies in its complementary set $\Omega^c$. \\

 Because $\Gamma$ is amenable, by F\o lner's criterion, one has 
 \begin{equation}
\label{eq-dem isop groupeamenable}
 	\infi{ \Omega} \ \frac{ \sharp \partial \Omega}{ \sharp \Omega} = 0\ 
\end{equation}
where the infimum ranges over all finite subsets $\Omega$ of $\Gamma$. \\

Being given a subset $\Omega \subset \Gamma$ we construct the analogous in $M$ as follows: fix $x_0 \in M$ any base point and let
	$$\Omega_M :=  \underset{\gamma \in \Omega}{\sqcup} D(\gamma \cdot x_0) \ . $$

Note that	
\begin{equation}
\label{eq-dem isop groupeamenable2}
	\mu_g(\Omega) = \vol(M_0) \cdot \sharp \Omega \ .
\end{equation}

To show that $\mathcal{O}$ does not have the spectral gap property we show that there is a constant $C > 0$ such that	for any $\Omega \subset \Gamma$ we have
	$$ \frac{ \left< (\Id - \mathcal{O})(\mathds{1}_{\Omega_M}) \ \cdot \mathds{1}_{\Omega_M})  \right> }{\sharp  \Omega_M } \le C \ \frac{ \sharp \partial \Omega}{ \sharp \Omega} \ , $$
which will conclude by using \eqref{eq-dem isop groupeamenable}. \\

We denote by $R$ the diameter of the manifold $M_0$. 
For a subset $U \subset M$ and $c > 0$ we denote by $N_{c}(U)$ the $c$-neighbourhood of $U$. \\

It is easy to verify that for all subsets $U \subset M$ one has 
\begin{itemize}
	\item  for any $x \in N_{R}(\partial U)$
 $$ |(\Id - \mathcal{O})(\mathds{1}_{U})(x)| \le 1 \ .$$ 
 	\item for any $x \notin  N_{R}(\partial U)$
 	 $$ (\Id - \mathcal{O})(\mathds{1}_{U})(x) = 0 \ .$$ 
 \end{itemize}

This readily yields constants $C_2, C_3 > 0$ such that  
\begin{align*}
	 \left< (\Id - \mathcal{O})(\mathds{1}_{\Omega_M}) \ \cdot \mathds{1}_{\Omega_M})  \right> & \le C_2 \cdot N_{R}(\partial \Omega_M) \\
	 	& \le C_3 \cdot \sharp( \partial \Omega)) \ , 
 \end{align*}
which, combined with \eqref{eq-dem isop groupeamenable2}, concludes. \hfill $\blacksquare$

\subsection{Reduction to the non-local isoperimetric inequality.} \label{subsec-reduction to isop} 
This  subsection is devoted to reducing the proof of Theorem \ref{theo-spectral gap} to the one of Proposition \ref{prop-isoperimetric inequality}. Recall that we want to prove that there is a constant $ c < 1$ such that for any $f \in L^2(M)$ we have 
	$$ || \mathcal{O}(f) ||_2 \le c ||f||_2 \ . $$
Since $\mathcal{O}$ is self-adjoint, the above inequality is equivalent to 
		$$ \supr{ f  \in L^2(M)} \ \frac{\left<\mathcal{O}^2(f) \cdot f \right>}{|| f||^2_{2}} < 1 \ , $$
which is what we will focus on proving from now on. \\

The proof is in two steps. The first one is to relate the spectral gap property of our operator to the energy associated to some quadratic form. Recall that we denoted by $p^{*2}(x,y)$ the kernel of the operator $\mathcal{O}^2$. Define  
$$D(f,f):= \frac{1}{2} \inte{M \times M} (f(y)-f(x))^2 \ p^{*2}(x,y) \ d \mu_g(x) d \mu_g(y)  \ .$$

\begin{lemma}
	For any smooth and compactly supported function $f : M \to \RR$ we have
\begin{equation}  
	\label{eq-demreduction prop isop 1}
		D(f,f) = \langle (\Id- \mathcal{O}^2) f, f\rangle  \ , 
\end{equation}
\end{lemma}

\textbf{Proof.} It follows from a simple computation which relies on the fact that $p^{*2}(x,y)$ is the kernel of a Markov operator.
\begin{eqnarray*}
\langle (I-\mathcal{O}^2) f, f\rangle &=& \inte{M} (I-\mathcal{O}^2)(f) \cdot f \ d \mu_g \\
 									&=& \inte{M} \left( f^2(x)-f(x)\inte{M} f(y) \ p^{*2}(x,y) \ d \mu_g(y) \right) \ d \mu_g(x) \\
 									&=& \inte{M}  f^2 \ d \mu_g  -\inte{M \times M} f(x)f(y)\ p^{*2}(x,y) \ d \mu_g(y)  d \mu_g(x) \\
 									&=& \frac{1}{2} \left(\inte{M}  f^2 \ d \mu_g -2 \inte{M \times M} f(x)f(y)\ p^{*2}(x,y) \ d \mu_g(y) d \mu_g(x) + \inte{M}  f^2  \ d \mu_g \right)\\
 									&=& \frac{1}{2} \left( \inte{M \times M} ( f^2(x)-2 f(x)f(y) + f^2(y)) \ p^{*2}(x,y) \ d \mu_g(x) d \mu_g(y) \right)\\
 									&=& D(f,f) \ .
\end{eqnarray*}
	
	\hfill $\blacksquare$ \\

Since $\mathcal{O}^2$ is self-adjoint and positive, we have 
	$$\| \mathcal{O}^2 \| = \sup_{f\in L^2(M)} \frac{\langle \mathcal{O}^2(f) , f\rangle }{\|f\|_2^2}=1-\infi{f\in L^2(M)}\frac{D(f,f) }{\|f\|_2^2} \ . $$

Therefore, Theorem \ref{theo-spectral gap} follows from the 
\begin{proposition}
	\label{prop.isopecheeger}
There exists $\epsilon>0$ such that for any $f \in L^2(M)$ we have
$$ \frac{D(f,f)}{  \|f\|^2_2} \geq \epsilon \ . $$
\end{proposition}	

The proof of the above Proposition  follows the same lines as Cheeger's proof of Cheeger's inequality. We chose to briefly recall its proof here since it is simpler and enlightens the reading of the proof in our setting. The reader familiar with it should perhaps skip to what is next. \\

\subsection{Interlude: Cheeger's inequality} The analogous of the non local isoperimetric inequality \ref{prop-isoperimetric inequality} in Cheeger's setting is simply given by the more classical isoperimetric inequality: for a non compact manifold $M$ we define its Cheeger's constant as 
			$$ h_1 := \inf_{\Omega} \ \frac{\mu_g^{d-1}(\partial \Omega)}{\mu_g(\Omega)} \ ,$$
where $\mu_g^{d-1}$ is the $d-1$ Haussdorf measure associated to the Riemannian metric $g$ and where the infinum is taken over all bounded subsets $\Omega$ of $M$ with smooth boundary. 

\begin{theoreme}[Cheeger]
	\label{theo.cheeger}
		For any smooth non zero compactly supported function $f : M \to \RR$ one has 
			$$ \frac{(\Delta f \cdot f)}{||f||_2^2} \ge \frac{h_1^2}{4} \ .$$ 
\end{theoreme}

\textbf{Proof.} By Stoke's formula, for any compactly supported function $f$ we have 
	\begin{equation}
		\label{eq.laplacestokes}
			 \left<\Delta f, f\right> = || \nabla f||^2_{2} \ .
	\end{equation}
	
We shall give the following pivotal quantity an upper bound and a lower one
	$$ I := \inte{M} |\nabla (f^2) | \ d \mu_g $$
By using the chain rule we get
	$$ I := \inte{M} 2 \ |\nabla f| \ |f| \ d\mu_g \ ,$$
which gives the upper bound by Cauchy-Schwartz inequality:
\begin{equation}
	\label{eq.cheegerupper}
	 I \le 2 \ || \nabla f||_{2} \ ||f||_{2}  \ . 
\end{equation}

The lower one is more subtle and starts with the use of the co-area formula:
\begin{align*}
	 I & = \inte{M} |\nabla (f^2) | \ d \mu_g \\
	 	& =  \inte{\RR_+} \mu^{n-1}_g(\{ f^2=t \})\ dt \ .
	\end{align*}
By construction of $h_1$ we get 
\begin{equation*}
		 	\inte{\RR_+} \mu^{d-1}_g(\{ f^2=t \})\ dt \  \ge h_1 \ \inte{\RR_+}   \mu_g(\{ f^2 \ge t \} ) \ dt = h_1 \ ||f||^2_{2} \ .
\end{equation*}

Combined with \eqref{eq.cheegerupper} we get
	$$ 2 \ || \nabla f||_{2} ||f||_{2} \ge h_1 \ ||f||^2_{2} \ .$$
and 
	$$ \frac{|| \nabla f||_{2}}{||f||_{2}} \ge \frac{h_1}{2} \ ,$$
which concludes by squaring both sides of the above equation and by using \eqref{eq.laplacestokes}. \hfill $\blacksquare$ \\

\textbf{Proof of (Proposition \ref{prop-isoperimetric inequality} $\Rightarrow$ Proposition \ref{prop.isopecheeger}).} As already emphasised, we shall mimic Cheeger's proof in our setting. There is essentially one point to handle: our operator does not come from the quadratic form $f \mapsto ||\nabla f||_2^2$ which prevents one to use the co-area formula. The trick to mimic the previous proof is to replace 
	$$ \inte{M} |\nabla f| \ d \mu_g$$
	 with 
$$	S(f) := \frac{1}{2} \inte{M \times M} |f(y)-f(x)| \ p^{*2}(x,y) \ d \mu_g(x) d \mu_g(y)  \ .$$ 

 The following lemma has to be compared with the end of the previous proof (from Equation \eqref{eq.cheegerupper} to the end) which addresses the use of the non local co-area formula together with Proposition \ref{prop-isoperimetric inequality}.
\begin{lemma}\label{lemme-cheeger}
Under the conclusion of Proposition \ref{prop-isoperimetric inequality}, there exists $\alpha>0$ such that for any $f \in L^2(M)$ we have
$$S(f^2) \geq \alpha \| f\|^2_2.$$
\end{lemma}

\textbf{Proof.} Using the symmetry in $x$ and $y$ we have
\begin{eqnarray*}
S(f^2) &=& \inte{ \{ (x,y) \in M^2 \ ,\  f^2(x)>f^2(y) \} } (f^2(x) -f^2(y)) \ p^{*2}(x,y) \ d \mu_g(x) d \mu_g(y) \\
 		&=& \inte{\RR} \inte{M \times M} \mathds{1}_{\{t, f^2(x)<t<f^2(y)\} } (t) \ p^{*2}(x,y) \ d \mu_g(x) d \mu_g(y)  \ dt \\
 		& = & \inte{\RR} \ \inte{U_t} \inte{U_t^c} \ p^{*2}(x,y) \ d \mu_g(x) d \mu_g(y)  dt \ ,
\end{eqnarray*}
where $U_t := \{ f^2 >  t\} $. \\

Now we apply the non local isoperimetric inequality given by Proposition \ref{prop-isoperimetric inequality}  to the set $U_t$ to get 
$$ S(f^2)  \geq  \inte{\RR} \beta \ \mu_g(U_t) \ dt = \beta \ \|f^2\|_1 = \beta \ \|f\|^2_2 \ , $$
concluding. \hfill $\blacksquare$ \\

Let us now adapt the first step in Cheeger's proof (relating $D(f,f)$ with $S(f^2) $) and conclude by using Lemma \ref{lemme-cheeger}.  \\

 Applying successively Cauchy-Schwarz inequality and the classical $(a+b)^2\leq 2(a^2+b^2)$ we get: 
\begin{eqnarray*}
S(f^2)^2  &=& \frac{1}{4} \left( \inte{M \times M}  |f(x) -f(y)| \ \sqrt{p^{*2}(x,y)}  \  |f(x) + f(y)| \  \sqrt{p^{*2}(x,y)} \ d \mu_g(x) d \mu_g(y) \right)^2 \\
			&\leq &\frac{1}{4} \inte{M \times M} (f(x)-f(y))^2 p^{*2}(x,y) \ d \mu_g(x) d \mu_g(y)     \inte{M\times M} (f(x)+f(y))^2 p^{*2}(x,y) \ d \mu_g(x) d \mu_g(y)  \\
			&\leq &\frac{D(f,f)}{2}  \inte{M \times M}( f(x)^2+f(y)^2) \ p^{*2}(x,y) \ d \mu_g(x) d \mu_g(y)  \ .     \\ 
\end{eqnarray*}
Using that $p(x,y)$ is the kernel of a Markov operator we also have
$$ \inte{M \times M}( f(x)^2+f(y)^2) \ p^{*2}(x,y) \ d \mu_g(x) d \mu_g(y)  = 2 \ \| f \|_{2}^2 \ , $$

and then 
\begin{equation} 
	\label{eq-demreduction prop isop 2}
S(f^2)^2 \leq D(f,f) \ \|f\|^2_2 \ .
\end{equation}

 Applying the conclusion of Lemma \ref{lemme-cheeger} yields
$$ S(f^2) \geq \alpha \|f^2\|_1  = \alpha \ \|f\|_2^2 \ , $$ 
which by using Inequality \eqref{eq-demreduction prop isop 2} and squaring the above inequality gives
 $$ \alpha^2 \ \|f\|_2^4 \le  S(f^2)^2 \leq D(f,f) \ \|f\|^2_2 \ .$$
Therefore,
 $$ D(f,f) \ge \alpha^2 \ \|f\|_2^2 \ ,$$
which is the desired inequality. \hfill $\blacksquare$

\section{Non-local isoperimetric inequality} 
\label{sec-proof of isop}

This section is devoted to the proof of Proposition \ref{prop-isoperimetric inequality}. We recall it here for the reader's convenience.

\begin{proposition}
	\label{prop-isoper dem}
There is $\beta>0$ such that for all relatively compact open sets $U\subset X$
$$ \inte{U} \inte{U^c} p^{*2}(x,y) \ d \mu_g(y) d \mu_g(x) \geq \beta \mu_g (U) \ .$$ \end{proposition}

The proof of the above proposition encounters two  difficulties of related but different nature. The point is that Proposition \ref{prop-isoper dem} does not reduce to an isoperimetric inequality as there is no constant $C > 0$ such that 	
\begin{equation}
	\label{eq-pas isoperi}
	 \inte{U} \inte{U^c} p^{*2}(x,y) \ d \mu_g(y) d \mu_g(x) \geq C \mu^{d-1}_g (\partial U) \ ,
\end{equation}
where $\mu^{d-1}_g$ stands for the $d-1$ Haussdorf measure.  Consider for example a very long and thin open set $U$ rolling around a disk of radius $1$ (Buser's hair). Dealing with 'hairy sets' was already the point of \cite{artbuser} along the proof of Buser's Inequality. One can also consider a Koch snowflake $U$: its boundary has infinite $d-1$ Haussdorf measure but the above left integral is bounded. The latter counter example to \eqref{eq-pas isoperi} is more about the non locality rather than the problem of Buser's hair and is something that should be dealt with during the proof. \\

However, our proof of Proposition \ref{prop-isoper dem} relies in the end on an isoperimetric inequality. The two following subsections aim at defining the main notion we will use to do so.

\subsection{Self-fat sets.} For the two next subsections we will assume that $M_0$ has injectivity radius greater or equal than $2$. We will often work at scale $1$, which is why \textbf{we introduce the notation $B(x) := B(x,1)$}. \\

		We say that a subset $A \subset M$ is $\alpha$-\textbf{self-fat} if for any $x \in A$ one has 
	$$ \mu_g(A \cap B(x)) \ge \alpha \ \mu_g(B(x)) \ . \\ $$ 

Note that the notion of self-fatness is highly dependent on $\alpha$. For example, take $\RR^2$ with the euclidean metric: for $\alpha \ge 1/2$ no balls of $\RR^2$ are self-fat regardless of the radius (because of the boundary points). If $\alpha< 1/2$ any sufficiently large ball becomes self-fat. \\

Note that the self-fatness is stable under union: if $A$ and $B$ are $\alpha$-self-fat then $A \cup B$ is also $\alpha$-self-fat. \\

Let us introduce some notations. For any $1\geq r > 0$ we define 
\begin{itemize}
	\item $v_-(r) := \infi{ x \in M} \ \mu_g(B(x,r)) \ ; $
	\item $v_+(r) := \supr{ x \in M} \  \mu_g(B(x,r)) \ . $ \\
\end{itemize}
Since we assumed that the injectivity radius of $M_0$ is greater than $1$ and because $\pi$ is a Riemannian covering we have
$$v_{-}(r) = \infi{ x \in M_0} \mu_g(B(x,r)) \ ,$$
and the analogous equality for $v_+(r)$. As a consequence, both the above defined functions of $r$ are continuous and positive (by compactness of $M_0$). 

\begin{definition}
	\label{def-selffat}
We define $$\alpha_0 := \frac{v_-(1/2)}{2 v_+(1)} \ $$ as the \textbf{fat parameter}. We say that a set $A$ is \textbf{self-flat} if it is $\alpha_0$-self-fat.
\end{definition}

\begin{remark} 
	\label{remark-selffat} This choice of $\alpha$ may seem arbitrary at this point. We shall see in the next subsection why we set it this way. Roughly, setting $\alpha_0$ as above, guarantees that the self-fat part (defined in the next subsection) of large balls is non empty. 
\end{remark}

If $A \subset M$ and $\epsilon > 0$ we denote by $A_{+, \epsilon}$ the '$\epsilon$-out neighbourhood of $A$' defined as
	$$ A_{+, \epsilon} := \{ \ x \in A^c \ , \ d(x,A) \le \epsilon \ \} \ .$$

The following lemma is one of the key to prove Proposition \ref{prop-isoper dem}. 

\begin{lemma}
\label{lemma-selffatset}
There are constants $C, \epsilon > 0$ such that for any self-fat set $A \subset M$ we have
	$$  \inte{A} \inte{A_{+,\epsilon}} \mathds{1}_{B(x)}(y) \ d \mu_g(x) d \mu_g(y) \ \ge C \ \mu_g(A_{+,\epsilon}) \ .$$
\end{lemma}

\begin{figure}[h!]
\begin{center}
	\def\svgwidth{0.4 \columnwidth}
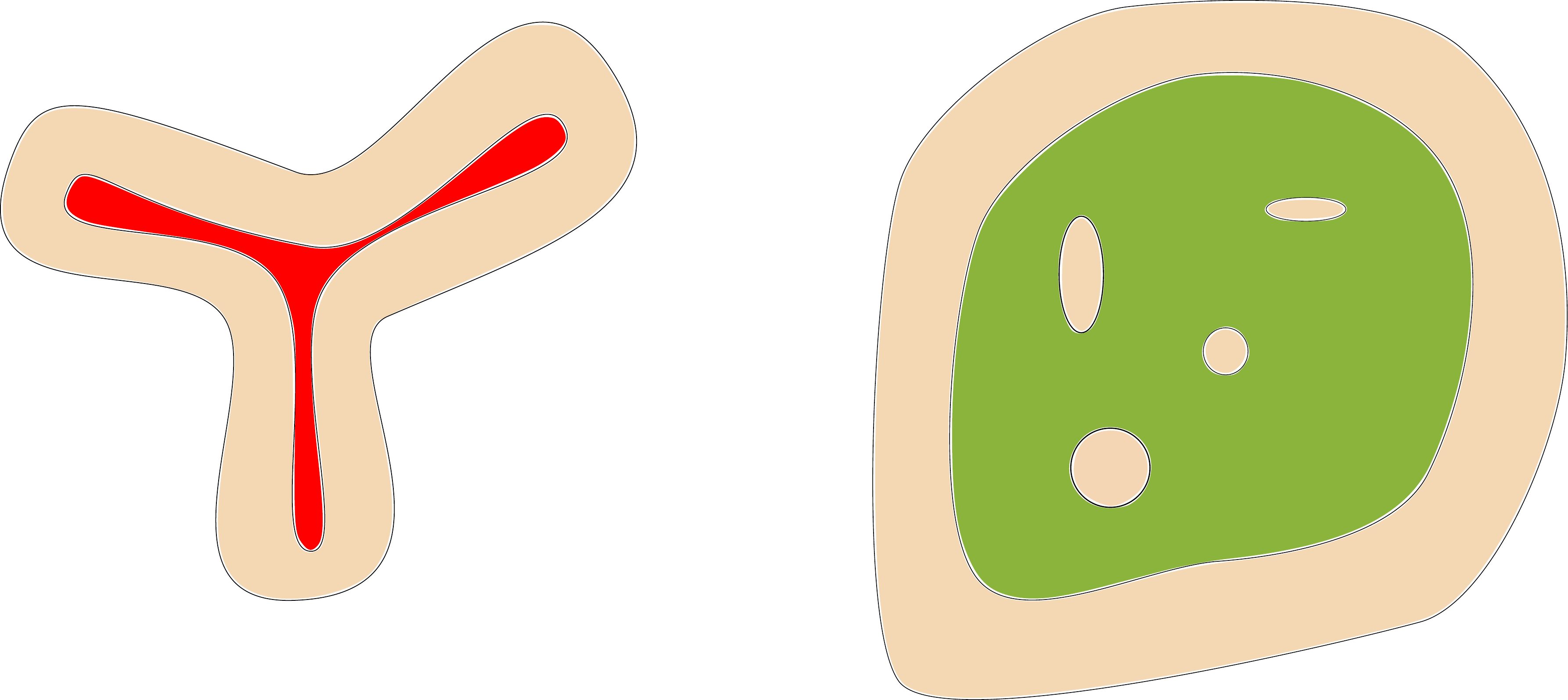
	\end{center}
\caption{The notion of self fat sets prevents the left picture to happen : too thin sets see their positive $\epsilon$-neighbourhood having more mass than themselves. On the contrary, we want to allow sets like the green one on the right.}
	\label{fig.selffat}
\end{figure}

\textbf{Proof.} We shall first adjust $\epsilon$. Consider the continuous function 
$$ \fonctionbis{[0,1]\times M }{\RR_+}{(\epsilon,x)}{  \supr{y \in B(x,\epsilon)} \mu_g(B(x) \Delta B(y))} \ \, $$
where $ A	\Delta B$ stands for the symmetric difference of the sets $A$ and $B$.
We Remark that this function takes the value $0$ when $\epsilon=0$, and is $\Gamma$-invariant with respect to its second variable. By compactness of $M/\Gamma$ there exists $\epsilon>0$ such that 
$$ \supr{y \in B(x,\epsilon)}   \mu_g(B(x) \Delta B(y)) \le \frac{\alpha_0 v_-(1)}{2} \ . $$

Let $A$ be a subset of $M$, let $x,y\in M$, one has:
$$\mu(B(x)\cap A) = \mu((B(x)\setminus B(y)) \cap A ) + \mu (B(x)\cap B(y) \cap A)$$
and
$$\mu(B(y)\cap A) = \mu((B(y)\setminus B(x)) \cap A ) + \mu (B(x)\cap B(y) \cap A).$$
Therefore, if $y\in B(x, \epsilon)$ one has : 
\begin{eqnarray*}
 \mu (B(y) \cap A) &=& \mu(B(x) \cap A) -\mu (B(x) \Delta B(y) \cap A) \\
 							&\geq & \mu(B(x) \cap A)-  \frac{\alpha_0 \ v_-(1)}{2}.		
\end{eqnarray*}
In particular for any self-fat set $\cA$ one has 

$$ \mu(B(x) \cap A) \ge  \alpha_0 \mu(B(x)) \ge \alpha_0 v_-(1) \ .$$
Therefore, for any  $y \in \cA_{+,\epsilon}$:
	$$ \mu_g(B(y) \cap \mathcal{A}) \ge \frac{\alpha_0 \ v_-(1)}{2} \ . $$

Finally we have: 
\begin{align*}
	\inte{\mathcal{A}} \inte{\mathcal{A}_{+, \epsilon}} \mathds{1}_{B(x)}(y) \ d \mu_g(y) d \mu_g(x) & =  \inte{\mathcal{A}_{+, \epsilon}}	\inte{\mathcal{A}} \mathds{1}_{B(y)}(x) \ d \mu_g(x) d \mu_g(y) \\
		& = \inte{\mathcal{A}_{+, \epsilon}} \mu_g(B(y) \cap \mathcal{A}) \ d \mu_g(y) \\
	& \ge \inte{\mathcal{A}_{+, \epsilon}} \frac{\alpha_0 \ v_-(1)}{2} \ d \mu_g(y)  \\	
		& \ge  \frac{\alpha_0 \ v_-(1)}{2} \cdot \mu_g(\mathcal{A}_{+, \epsilon}) \ ,
\end{align*}
which concludes the proof. \hfill $\blacksquare$

\subsection{Self-fat part of sets.} This subsection aims at investigating the following notion.

\begin{definition}
\label{def-selffatpart}
Given a set $U \subset M$ we define its \textbf{self-fat part}, which we denote $\mathcal{SF}(U)$, as the maximal self-fat set contained in $U$. It is well defined since, as previously noticed, the self-fatness is stable under union. 
\end{definition}

The second key property we will need is contained in the following lemma. In order to state it, let us introduce another kind of subsets of a given set $U$. \\ 

For any number $ \kappa \in [0,1]$ we call the \textbf{$\kappa$-thick part} of $U$, that we denote $\mathcal{E}_{\kappa}$, the set 
	$$ \mathcal{E}_{\kappa}(U) := \{ \ x \in U \ , \ \mu_g(U \cap B(x)) \ge \kappa \cdot \mu_g(B(x)) \ \} \ . $$

\begin{lemma}
	\label{lemma-selffatsubset}
There is $0 < \kappa < 1$ such that for any $U \subset M$ we have 
	$$ \mathcal{E}_{\kappa}(U) \subset \mathcal{SF}(U)\ .$$
\end{lemma}

Following up with Remark \ref{remark-selffat}, the above lemma implies in particular that 
$$B(x,r-1) \subset \mathcal{SF}(B(x,r))$$ 
for any $r > 1$. This justifies the choice of $\alpha_0$. \\

\textbf{Proof.} The constant $\kappa > 0$ is chosen such that for any $x \in M$ and any $U \subset M$ if 
\begin{equation}
	\label{eqref-selffatpart1}
	 \mu_g(B(x) \cap U) \ge \kappa  \cdot \mu_g(B(x)) 
\end{equation}
then 
\begin{equation}
	\label{eqref-selffatpart2}
		\mu_g(B(x, 1/2) \cap U) \ge \frac{\mu_g(B(x, 1/2))}{2} \ . 
\end{equation}

Let us see why such a $\kappa$ exists. Indeed, on the one hand we have 
$$ \mu_g(B(x, 1/2) \cap U) + \mu_g(B(x, 1/2)^c \cap B(x) \cap U) = \mu_g(B(x) \cap U) \ge \kappa  \cdot \mu_g(B(x)) \ ,$$
under Inequality \eqref{eqref-selffatpart1}. On the other hand 
\begin{align*}
	\mu_g(B(x, 1/2)^c \cap B(x) \cap U) & \le \mu_g(B(x, 1/2)^c \cap B(x)) \\
		& \le \mu_g(B(x)) - \mu_g(B(x,1/2)) \ .
\end{align*}

Therefore, 
$$ \mu_g(B(x, 1/2) \cap U) + \mu_g(B(x)) - \mu_g(B(x,1/2)) \ge \kappa  \cdot \mu_g(B(x)) \ ,$$
which can be rewritten as 
$$ \mu_g(B(x, 1/2) \cap U) \ge \mu_g(B(x,1/2)) + (\kappa -1) \cdot \mu_g(B(x)). $$

Setting $\kappa$ close enough to one in order for the following to hold
	$$ (1- \kappa) v_+(1) \le \frac{v_-(1/2)}{2} , $$
we get 
	$$ \mu_g(B(x, 1/2) \cap U) \ge \frac{\mu_g(B(x,1/2))}{2} \ .$$

Let us now see how such a choice of $\kappa$ implies that the conclusion of Lemma \ref{lemma-selffatsubset} holds. \\

Let $U$ be any subset of $M$ and $x \in \mathcal{E}_{\kappa}(U)$. We want to show that $x \in \mathcal{SF}(U)$. We will actually show that $B(x, 1/2) \cap U \subset \mathcal{SF}(U)$. Since $\mathcal{SF}(U)$ is maximal for self-fatness, we want to show that $B(x, 1/2) \cap U$ is self-fat. \\

In other word, we want to show that for any $y \in B(x, 1/2) \cap U$ 
	$$ \mu_g(B(y) \cap B(x, 1/2) \cap U) \ge \alpha_0 \ \mu_g(B(y)) \ .$$

Note that $B(x, 1/2) \subset B(y)$ for any $y \in B(x, 1/2)$. In particular for any $y \in B(x, 1/2) \cap U$ we have
	$$  \mu_g(B(y) \cap B(x, 1/2) \cap U) =  \mu_g( B(x, 1/2) \cap U) \ .  $$ 
Because we supposed that $x \in \mathcal{E}_{\kappa}(U)$ and because of our choice of $\kappa$ we have by Inequality \eqref{eqref-selffatpart2}
\begin{align*}
	\mu_g( B(x, 1/2) \cap U)  & \ge  \frac{\mu_g(B(x, 1/2))}{2} \\
			& \ge \frac{v_-(1/2)}{2} \ .
\end{align*}
We conclude by using our choice of self-fat parameter $\alpha_0$:
\begin{align*}
\mu_g( B(x, 1/2) \cap U) &  \ge  \alpha_0 \ v_+(1) \\
		&  \ge  \alpha_0 \ \mu_g(B(y)) \ , 
\end{align*}
by definition of $v_+(1)$.  \hfill $\blacksquare$

\subsection{Proof of Proposition \ref{prop-isoper dem}} We conclude this section by the proof of Proposition \ref{prop-isoper dem}. \\

Let $\epsilon_0 := \infi{ x \in M} \mathrm{inj}_M(x)$ where $\mathrm{inj}_M(x)$ is the injectivity radius at $x$ of the manifold $M$. Note that $\epsilon_0 > 0$ since we supposed $M_0$ compact and because $ \epsilon_0 \ge \infi{ x \in M_0} \mathrm{inj}_{M_0}(x)$. Note also that Proposition \ref{prop-isoper dem} is invariant under metric scaling; given $\lambda > 0$ Proposition \ref{prop-isoper dem} holds for the metric $g$ if and only if it holds for the metric $\lambda^2 g$. Up to using such a scaling one can (and one will) suppose $\epsilon_0 \ge 2$, making the statements of the last subsection to fit in. \\

It follows from the construction of the Dirichlet domain at $x$ that
	$$ p(x,y) \ge \frac{1}{\mu(M_0)} \mathds{1}_{B(x,1)}(y) =  \frac{1}{\mu(M_0)} \mathds{1}_{B(x)}(y) \ ,$$
	
since we supposed the injectivity radius to be greater than $2$. Note that by construction one has 
$$ p^{*2}(x,y) = \inte{M} \ p(x,z) \ p(z,y) \ d \mu_g(z) \ . $$ 
In particular, there is a constant $c > 0$ such that for all $x,y \in M$ 
$$ p^{*2}(x,y) \ge c \ \mathds{1}_{B(x)}(y) \ .$$

Therefore, Proposition \ref{prop-isoper dem} follows from the following statement. \\

\textit{There is $\beta>0$ such that for all relatively compact open sets $U\subset M$}
\begin{equation}
	\label{eqref-propdem isoper}
	 \inte{U} \inte{U^c} \mathds{1}_{B(x)}(y) \ d \mu_g(y) d \mu_g(x) \geq \beta\ \mu_g(U) \ . 
\end{equation} 

We will now focus on proving that the above inequality holds. In order to do so, we first split $U$ into its self-fat part and its complement. \\

 We split the left member of Equation \eqref{eqref-propdem isoper} as \\
\begin{equation}
	\label{eq1-propdem isoper}
	\begin{split}
		 \inte{U} \inte{U^c} \mathds{1}_{B(x)}(y) \ d \mu_g(y) d \mu_g(x) = & \\
			   \inte{\mathcal{SF}(U)} \ \inte{U^c} \mathds{1}_{B(x)}(y) \ d \mu_g(y) d \mu_g(x) & +  \inte{U \setminus \mathcal{SF}(U)} \ \inte{U^c} \mathds{1}_{B(x)}(y) \ d \mu_g(y) d \mu_g(x) \ . 
	\end{split} 
\end{equation} \\

We will bound from below the two above integrals independently. Let us start by bounding from below the last one, which is the easiest to deal with. \\

Recall the conclusion of Lemma \ref{lemma-selffatsubset} which asserts that there is $ 0<\kappa < 1$ such that
	$$ \mathcal{E}_{\kappa}(U) \subset \mathcal{SF}(U)  \ . $$ 

Taking the complementary set (as subsets of $U$) we get
	$$ U \setminus \mathcal{SF}(U) \subset U \setminus \mathcal{E}_{\kappa}(U)  \ . $$ 

Because of how $\mathcal{E}_{\kappa}(U)$ is defined we have
\begin{align*}
	 U \setminus \mathcal{E}_{\kappa}(U) & = \{ \ x \in U \ , \ \mu_g(U \cap B(x)) \le \kappa \cdot \mu_g(B(x)) \ \} \\ 
	 & =  \{ \ x \in U \ , \ \mu_g(U^c \cap B(x)) \ge (1- \kappa) \cdot \mu_g(B(x)) \ \} 
\end{align*}

In other word, from the perspective of a point $x \in  U \setminus \mathcal{SF}(U)  \subset U \setminus \mathcal{E}_{\kappa}(U)
$, some definite mass of its $1$-neighbourhood is carried by $U^c$. In particular, for any $x \in U \setminus \mathcal{SF}(U) $ one has 
\begin{align*}
	\inte{U^c} \mathds{1}_{B(x)}(y) \ d \mu_g(y) & = \mu_g(B(x)\cap U^c) \\
		& \ge (1- \kappa) \cdot \mu_g(B(x)) \\
		& \ge (1- \kappa) \cdot v_-(1) \ .
\end{align*}

Therefore,
\begin{equation}
	\label{eq3-propdem isoper}
	   \inte{U \setminus \mathcal{SF}(U)} \ \inte{U^c} \mathds{1}_{B(x)}(y) \ d \mu_g(y) d \mu_g(x) \ge (1- \kappa) v_-(1) \cdot \mu_g(U \setminus \mathcal{SF}(U))  \ .
\end{equation}

Let us now bound from below the other integral appearing in the bottom of Equation \eqref{eq1-propdem isoper}. \\

We fix from now on a pair $C, \epsilon > 0$ satisfying the conclusion of lemma \ref{lemma-selffatset}. We start with the two obvious lower bounds valid for any $ 0 \le \delta \le 1$

\begin{align*}
	 \inte{\mathcal{SF}(U)} \ \inte{U^c} \mathds{1}_{B(x)}(y)  \ d \mu_g(y) d \mu_g(x) & \ge \inte{ \mathcal{SF}(U)} \ \inte{\mathcal{SF}(U)_{+, \epsilon} \cap U^c} \mathds{1}_{B(x)}(y) \ d \mu_g(y) d \mu_g(x) \\
	 & \ge \delta  \inte{\mathcal{SF}(U)} \ \inte{ \mathcal{SF}(U)_{+, \epsilon} \cap U^c} \mathds{1}_{B(x)}(y) \ d \mu_g(y) d \mu_g(x) \ . \\
\end{align*}

Because of 
$$ \mathcal{SF}(U)_{+, \epsilon}  = (\mathcal{SF}(U)_{+, \epsilon} \cap U) \sqcup (\mathcal{SF}(U)_{+, \epsilon} \cap U^c) $$ we have for any $\delta \in [0,1]$: 
\begin{equation}
	\label{eq2-propdem isoper}
		\begin{split}
 \inte{\mathcal{SF}(U)} \ \inte{U^c} \mathds{1}_{B(x)}(y)  \ d \mu_g(y) d \mu_g(x)  & \ge \\ 
 \delta \Big( \inte{\mathcal{SF}(U)} \ \inte{\mathcal{SF}(U)_{+, \epsilon}} \mathds{1}_{B(x)}(y)  \ d \mu_g(y) d \mu_g(x) & -  \inte{\mathcal{SF}(U)} \ \inte{ \mathcal{SF}(U)_{+, \epsilon} \cap U} \mathds{1}_{B(x)}(y)  \ d \mu_g(y) d \mu_g(x) \Big) \ .
	 	 \end{split}
\end{equation} \\

We shall first give a lower bound to the last above integral using the following rough upper bounds
\begin{align*}
	 \inte{\mathcal{SF}(U)} \ \inte{\mathcal{SF}(U)_{+, \epsilon} \cap U} \mathds{1}_{B(x)}(y)  \ d \mu_g(y) d \mu_g(x) & = \inte{\mathcal{SF}(U)_{+, \epsilon} \ \cap U} \inte{\mathcal{SF}(U)}  \mathds{1}_{B(x)}(y)  \ d \mu_g(x) d \mu_g(y) \\
	  & \le \inte{\mathcal{SF}(U)_{+, \epsilon} \cap U} \mu_g(B(y))  \ d \mu_g(y)  \\
	  & \le \inte{y \in U \setminus \mathcal{SF}(U)} \mu_g(B(y))  \ d \mu_g(y) \\
	  & \le v_+(1) \cdot \mu_g(U \setminus \mathcal{SF}(U)) \ .
\end{align*}

Combining the above upper bound with \eqref{eq1-propdem isoper}, \eqref{eq3-propdem isoper} and \eqref{eq2-propdem isoper} we get for all $\delta \in [0, 1]$
\begin{equation*}
	\begin{split}
		 \inte{U} \inte{U^c} \mathds{1}_{B(x)}(y) & (d(x,y))  \ d \mu_g(y) d \mu_g(x) \ge (1- \kappa) v_-(1) \cdot \mu_g(U \setminus \mathcal{SF}(U)) \\ 
			 & +  \delta \Big( - v_+(1) \cdot \mu_g(U \setminus \mathcal{SF}(U)) + \inte{\mathcal{SF}(U)} \ \inte{\mathcal{SF}(U)_{+, \epsilon}} \mathds{1}_{B(x)}(y)  \ d \mu_g(y) d \mu_g(x) \Big) \ .
	\end{split} 
\end{equation*}

From now on, we fix $\delta \in [0,1]$ small enough such that there is $C_2 > 0$ such that  
	$$ (1- \kappa) v_-(1)  - \delta  v_+(1) > C_2 \ . $$
Which gives 
\begin{equation*}
	\begin{split}
		 \inte{U} \ \inte{U^c} \mathds{1}_{B(x)}(y) \ d \mu_g(y) d \mu_g(x) & \ge C_2 \cdot \mu_g(U \setminus \mathcal{SF}(U)) + \delta \inte{ \mathcal{SF}(U)} \ \inte{\mathcal{SF}(U)_{+, \epsilon}} \mathds{1}_{B(x)}(y) \ d \mu_g(y) d \mu_g(x) \ .
	\end{split} 
\end{equation*}

We conclude this proof by proving 
\begin{lemma}
	\label{lemme-propdem isoper}
	There is a constant $C_3 > 0$ such that for any open relatively compact set $U$ we have
		$$ \inte{\mathcal{SF}(U)} \ \inte{ \mathcal{SF}(U)_{+, \epsilon}} \mathds{1}_{B(x)}(y) \ d \mu_g(y) d \mu_g(x) \ge C_3 \ \mu_g(\mathcal{SF}(U)) \ .$$ 
\end{lemma}

The above lemma implies Proposition \ref{prop-isoper dem} by setting $\beta := \min \{ \delta C_3, C_2 \}$:
\begin{equation*}
	\begin{split}
		 \inte{U} \inte{ U^c} \mathds{1}_{B(x)}(y) \ d \mu_g(y) d \mu_g(x)  & \ge C_2 \mu_g(U \setminus \mathcal{SF}(U)) \ +  \delta C_3 \ \mu_g(\mathcal{SF}(U)) \\
		 & \ge \beta \ (\mu_g(U \setminus \mathcal{SF}(U)) + \mu_g(\mathcal{\mathcal{SF}(U)})) \\
		 & \ge \beta \ \mu_g(U) \ . 
	\end{split} 
\end{equation*}

\textbf{Proof of Lemma \ref{lemme-propdem isoper}.} Because of our choice of the pair $(\epsilon,C)$ and since is $\mathcal{SF}(U)$ is self-fat by construction  we have access to lemma \ref{lemma-selffatsubset} whose conclusion is 
$$ \inte{\mathcal{SF}(U)} \ \inte{\mathcal{SF}(U)_{+, \epsilon}} \mathds{1}_{B(x)}(y) \ d \mu_g(y) d \mu_g(x)  \ge C \ \mu_g(\mathcal{SF}(U)_{+, \epsilon}) \ .$$ 	

It remains then to prove that 	there is a constant $C_4$ such that 
$$ \mu_g(\mathcal{SF}(U)_{+, \epsilon}) \ge C_4 \cdot \mu_g(\mathcal{SF}(U)) \ . $$

We start off using the co-area formula with the function 
$$ \fonction{\Phi}{ \mathcal{SF}(U) \sqcup \mathcal{SF}(U)_{+, \epsilon}}{[0, \epsilon]}{x}{d(x, \mathcal{SF}(U))} \ , $$
which (as a distance of a given set) satisfies $|\nabla \Phi| =1$ almost everywhere on $\mathcal{SF}(U)_{+, \epsilon}$. Therefore, 
\begin{align*}
	 \mu_g(\mathcal{SF}(U)_{+, \epsilon}) & =  \inte{\mathcal{SF}(U)_{+, \epsilon}} |\nabla \Phi| \ d \mu_g \\
	 	& = \inte{]0,\epsilon[} \mu_g^{d-1}( \{ \Phi =t \} ) dt \ ,
\end{align*}
where $\mu_g^{d-1}$ is the codimension 1 Riemannian measure. \\ 

For what follows, we keep using the notations introduced in Subsection \ref{subsec-if and only if}. Since we assume that $\Gamma$ is non amenable, F\o lner's criterion asserts that 
\begin{equation}
\label{eq-dem isop groupe}
 	\infi{ \Omega} \ \frac{ \sharp \partial \Omega }{ \sharp \Omega} > 0 \ ,  
\end{equation}
where the infimum ranges over all subsets $\Omega$ of $\Gamma$.\\

In order to 'pull-back' \eqref{eq-dem isop groupe} to the manifold $M$ we shall use the following theorem which specifies in our setting as

\begin{theorem}{\cite[Corollaire 6.7]{artcoulhonsaloffisometrie}}
\label{theo-isopersaloff}
Let $M \to M_0$ be a Riemannian covering of deck group $\Gamma$ with $M_0$ compact without boundary. Then Inequality \eqref{eq-dem isop groupe} is equivalent to
	$$ \infi{ \Omega} \ \frac{ \mu_g^{d-1}(\partial \Omega)}{\mu_g(\Omega)} \ge C_5 > 0\ ,$$
where $\mu_g^{d-1}$ is the codimension one Haussdorf measure and where $\Omega$ ranges over open subsets of $M$ with regular boundary.
\end{theorem}

\begin{remark}
The authors of \cite{artcoulhonsaloffisometrie} do not precise their definition of 'regular'. One can consider equivalently sets of smooth boundary in the range of the above infimum or sets of finite perimeter since these latter are well approachable by these first. We refer to \cite[Chapter 5]{livretheoriegeommeasure} for more details on sets of finite perimeters. We want to use Theorem \ref{theo-isopersaloff} with level sets of a Lipschitz function which have (locally) almost surely finite codimension 1 Haussdorf measure. The perimeter is always less or equal than the Haussdorf measure (see \cite[Section 5.7]{livretheoriegeommeasure}) which allows us to use the above theorem. We recommend the first part of the \cite{surveynonlocalperimeter} for an introduction to the basics of geometric measure theory needed in this article.
\end{remark}

In particular for almost every $t \in [0, \epsilon]$ one has
\begin{align*}
	\mu_g^{d-1}( \{ \Phi =t \} ) & \ge C_5 \ \mu_g( \{ \Phi \le t \}) \\
		 & \ge C_5 \ \mu_g( \{ \Phi \le 0 \} ) \\
		 & \ge C_5 \ \mu_g( \mathcal{SF}(U)) \ .
\end{align*}
And then, 
\begin{align*}
	\mu_g(\mathcal{SF}(U)_{+, \epsilon}) & = \inte{]0,\epsilon]} \mu_g^{d-1}( \{ \Phi =t \} ) dt  \\ 
		& \ge C_5 \ \epsilon \ \mu_g( \mathcal{SF}(U)) \ , 
\end{align*}
concluding.  \hfill $\blacksquare$  $\blacksquare$

\section{Asymptotic behaviour of Dirichelt random walks}\label{sec-central_limit_theorem}

In this section we explain the geometric consequences of the spectral gap  and prove Theorems \ref{th-convergence dans le bord intro} and \ref{th-main-asymptotic-behaviour-of-Zn}. We shall see that Theorem \ref{th-convergence dans le bord intro} is a consequence of Kingman's subadditive ergodic theorem \cite{kingman} together with the so called linear progress with exponential tail property. 
Recall that  $M_0$ is  a compact manifold, endowed with the  probability measure given by the normalized Riemannian measure $\mu_g$.  Let $\Omega := M_0^\N$ and  $X_n$ be the coordinate maps from $\Omega$ to $M_0$ defined by $X_n(\omega) = \omega_n$ for a sequence $\omega=(\omega)_{i\in \N} \in \Omega$. We equip $\Omega$ with the product measure $\mu :=\mu_g^{\otimes \N}$. We denote by $T: \Omega \tv \Omega $ be the shift operator.\\

For $x \in M$ and $y \in M_0$, recall that we denoted by $D_x(y)$ the unique lift of $y$ in $M$ in the Dirichlet domain associated to $x$. Recall also that $Z_n$ is constructed recursively as $D_{Z_{n-1}}(X_n)$ with $Z_0=o$.

\subsection{Linear progress with exponential tail property.}
\label{subsec.linearprogresswithexp}
This section aims at showing that $Z_n$ satisfies the linear progress with exponential tail property that we recall here for the reader's convenience.

\begin{definition}
	A sequence $Z_n$ of random variables taking value in a metric space $(X,d)$ is said to satisfy the linear progress with exponential tail property if there is $C, \epsilon > 0$ such that for any $n,m \in \NN$ we have 
		$$ \PP( d( Z_n, Z_{n+m}) \le \epsilon m ) \le C \ e^{-\epsilon m} \ .$$
\end{definition}

The following proposition is a consequence of Theorem \ref{theo-spectral gap}.

\begin{proposition}
\label{prop-linearprogress}
Suppose the deck group $\G$ is non amenable then the sequence $Z_n$ has the linear progress with exponential tail property.
\end{proposition}

\textbf{Proof.} The proof proposed here follows the classical line of work as it relies on the spectral gap. We will however have to combine it with the following lemma that can be understood as an almost translation-in-time-invariance of the random variable $d(Z_n,Z_m)$. For the sake of clarity, we postpone the proof of this lemma to the end of this subsection. We will reuse this lemma all along this section.  Recall that $R$ is the diameter of $M_0$.

\begin{lemma}
\label{lemme.comparisondistance}
For any $n, m \in \NN$ we have 
	$$ | d(Z_n, Z_{n+m}) - d(o, Z_{m} \circ T^n) | \le R \ . $$
\end{lemma}

\textbf{Proof of (Lemma \ref{lemme.comparisondistance} $\Rightarrow$ Proposition \ref{prop-linearprogress}).} Let $0 < \epsilon_2 < \epsilon$ where $\epsilon$ is as in Theorem \ref{theo-spectral gap}. Lemma \ref{lemme.comparisondistance} gives that
	\begin{eqnarray*}
			\PP( d(Z_n, Z_{n+m}) \le \epsilon_2 m) & \le & \bP( d(o, Z_{m} \circ T^n) \le \epsilon_2 m + R) \\
&\le & \bP( d(o, Z_{m} \circ T^n) \le \epsilon m)  \ ,
	\end{eqnarray*}
for any $n$ and $m$ such that $m > C_1$ for some constant $C_1$ independent of $m$ and $n$. \\

Since $(X_i)$ are I.I.D random variables, $Z_m \circ T^n$ follows the law of $Z_m$. In particular,
	$$ \bP( d(o, Z_{m} \circ T^n) \le \epsilon m) =\bP( d(o, Z_{m}) \le \epsilon m) \ .$$

Proposition \ref{prop-linearprogress} then follows from the spectral gap theorem. Indeed, we show  that there is $C_2$ such that for all $\alpha > 0$ small enough we have that for any $m \ge 0$ 
\begin{equation}
	\label{eq.trouspectraluniform}
	 \PP(d(o, Z_{m}) \leq \alpha m) \le C_2 \ e^{- \alpha m} \ .
\end{equation}

Note first that by construction of the Dirichlet operator we have 
	\begin{eqnarray*}
		\PP(d(o, Z_{m}) \leq \alpha m) & = & \inte{B(o, \alpha n)} p^{*m}(o,y) \ d \mu_g(y) \\
		& = & \mathcal{O}^{n}(\mathds{1}_{B(o, \alpha n)})(o) \ .
	\end{eqnarray*}	 
	
To reach the Inequality \eqref{eq.trouspectraluniform} we will need some spatial uniformity that we get from the following simple remark:
for all $f\in L^2(M)$ and for all $x \in M$:
	\begin{equation}
		\label{eq.cauchyschwartzspatial}
		 |\cO(f)(x )| \leq \frac{1}{\sqrt{\Vol(M_0)}}\|f\|_2.
	\end{equation}
	
This is a simple consequence of Cauchy-Schwarz inequality: 
\begin{eqnarray*}
\cO (f)(x)& =&  \frac{1}{\mu_g(M_0)} \inte{M_0} 1_{D_x}(y) f(y) \ d \mu_g(y) \\
&\leq &  \frac{1}{\mu_g(M_0)} \sqrt{\inte{M_0} 1_{D_x}(y) \ d \mu_g(y)} \sqrt{ \inte{M_0} f(y)^2 \ d \mu_g(y)} \\
&\leq & \frac{1}{\sqrt{\mu_g(M_0)}}\|f\|_2.
\end{eqnarray*}

Let us come back to the proof of Inequality \eqref{eq.trouspectraluniform}. \\ 

Note first that $M$ has volume growth at most exponential since $M$ is roughly isometric, in the meaning of Kanaï, to some Cayley graph of $\Gamma$ (a finitely generated group has at most exponential growth). Let then $h > 0$ such that for any $x \in M$ and any $r > 0$ 
	$$ \mu_g(B(x,r)) \le h^{-1} e^{h r} \ .$$

Using Successively Inequality \eqref{eq.cauchyschwartzspatial}, Theorem \ref{theo-spectral gap} and the at most exponential volume growth we get 
	\begin{align*}
			 \PP(d(o, Z_{n}) \leq \alpha n) & = \cO^n(\mathds{1}_{B(o, \alpha n)})(x) \\
			& \leq  \frac{1}{\sqrt{\mu_g(M)}}\|\cO^{n-1}(\mathds{1}_{B(o, \alpha n)})\|_2 \\
			& \leq C_3(1-\epsilon)^n \|\mathds{1}_{B(o, \epsilon n)}\|_2  \\
			& \leq  C_3(1-\epsilon)^n \mu_g(B(x, \alpha n))  \\
			& \le  C_4 (1-\epsilon)^n e^{ \alpha h n}  \ ,
	\end{align*}
for some positive constants $C_3, C_4$. Choosing $\alpha > 0$ small enough such as 
	$$  - \alpha >  \ln(1 - \epsilon) +  \alpha h   \ $$
concludes the proof. \hfill $\blacksquare$ \\

\begin{figure}[h!]
\begin{center}
	\def\svgwidth{0.5 \columnwidth}
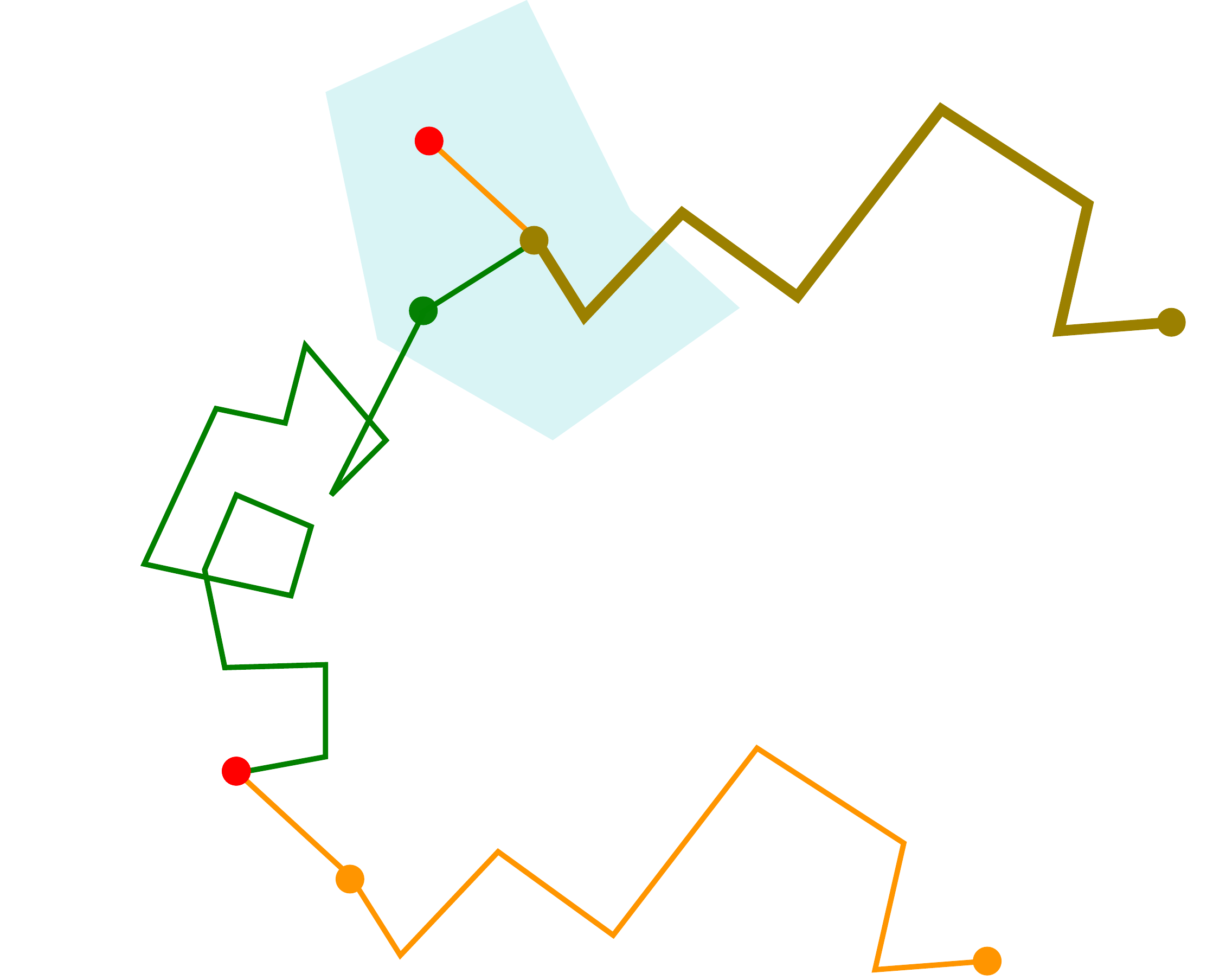
	\end{center}
\caption{}
	\label{fig.almostsubadditive}
\end{figure}

\textbf{Proof of Lemma \ref{lemme.comparisondistance}.} The proof is contained in Figure \ref{fig.almostsubadditive} that we comment bellow. \\

By construction, the Dirichlet domain $D_{Z_{n+1}}$, in blue in Figure \ref{fig.almostsubadditive}, contains $Z_n$ as well as some lift $\gamma \cdot o$ of $p_o$ for some $\gamma \in \Gamma$. The path relating $o$ to $Z_{m} \circ T^n$, in orange, is mapped isometrically by $\gamma^{-1}$ to the path relating $\gamma \cdot o$ to $Z_{m+n}$ in orange/brown. By construction of the Dirichlet domain both $Z_n$ and $\gamma \cdot o$ are $R/2$ close to $Z_{n+1}$ and therefore $Z_n$ and $\gamma \cdot 0$ are $R$ close from one another, concluding. 
\hfill $\blacksquare$ $\blacksquare$

\subsection{Proof of Theorem \ref{th-convergence dans le bord intro}.} This section is dedicated to the proof of Theorem \ref{th-convergence dans le bord intro} that we recall here for the reader's convenience.

\begin{theorem} \label{theroem convergence of d(o,Xn)/n}
Suppose the deck group $\G$ is non amenable then there exists $\ell>0$ such that for almost all sequence $Z_n$ one has: 
$$\lim_{n\tv \infty} \frac{d(o,Z_n)}{n}=\ell.$$
If moreover $\G$ is Gromov hyperbolic then the sequence $Z_n$ converges almost surely to a point in the Gromov boundary $\partial M$.
\end{theorem}

\textbf{Proof.} The part showing that $\ell > 0$ follows easily from Proposition \ref{prop-linearprogress}. Indeed, the Borel-Cantelli lemma implies that there is $\alpha > 0$ such that for almost all trajectories
 $$ \limiinf{n \tv \infty} \ \frac{d(o,Z_n)}{n}\geq c>0 \ . $$

A classical consequence of the above convergence is the second part of Theorem \ref{theroem convergence of d(o,Xn)/n}: almost all trajectories $(Z_n)_{n\in \N}$ converge in the Gromov boundary $\partial M$. Recall that $(Z_n)$ converges almost surely in $\partial M$ if and only if 
	\begin{equation}
		\label{eqdef.convergenceboundary}
		 \langle Z_n , Z_{m} \rangle_o  \tends{n,m \to \infty} \infty \  
	\end{equation}
almost surely. The proof that $Z_n$ converges almost surely in the boundary relies on the fact that (see \cite[7.3]{GdlH}) there exists an angle like function $\rho : M \times M \to \RR_+$ which is
	\begin{enumerate}
		\item compatible with the Gromov product: there is $\epsilon > 0$ and a constant $C > 0$ such that for all $x, y \in M$
		$$   C^{-1} \ e^{ - \epsilon \langle x , y \rangle_o}  \le  \rho(x,y) \le C \ e^{ - \epsilon \langle x , y \rangle_o} \ ; $$	
		\item a pseudo distance: for any $x,y,z \in M$ we have the triangular inequality 
			$$ \rho(x,z) \le \rho(x,y) + \rho(y,z) \ .$$
	\end{enumerate}
 
The second item above implies in particular that $\rho$ is a pseudo-distance. Even though we will not use it, this pseudo-distance extends as an actual distance to the Gromov boundary of $M$. \\

The proof then follows from Sullivan's one \cite{sullivanconvergenceboundary} and starts with an upper bound. We assume without any loss of generality that $m \ge n$ and we set $p := m -n$. Using both the above items we get
\begin{eqnarray}
	 e^{ - \epsilon \langle Z_n , Z_{m} \rangle_o} =  e^{ - \epsilon \langle Z_n , Z_{n+p} \rangle_o} & \le & C \ q(Z_n , Z_{n+p}) \\
	 & \le & C \somme{ n \le k \le n +p -1} q(Z_k , Z_{k+1}) \\
	 & \le & C^2 \somme{ n \le k \le n +p -1} e^{ - \epsilon \langle Z_k , Z_{k+1} \rangle_o} \\
	 & \le & C^2 \somme{ n \le k } e^{ - \epsilon \langle Z_k , Z_{k+1} \rangle_o} \ . 
\end{eqnarray}
Convergence \eqref{eqdef.convergenceboundary} then follows if one shows that 
	\begin{equation}
		\label{eq.convergencetotheboundnary}
		 \EE \left( \somme{n \in \NN} e^{ - \epsilon \langle Z_k , Z_{k+1} \rangle_o} \right) < \infty \ ,
	\end{equation}
since it implies that the series is almost surely finite. Since all the terms of the series are positive, it is sufficient to show that there is are constants $\alpha, C_2$ such that
	$$ \EE \left( e^{ - \epsilon \langle Z_k , Z_{k+1} \rangle_o} \right) \le C_2 \ e^{- \alpha k} \ .$$
This follows readily from the linear escape with exponential tail property since it implies in particular that there is $\alpha > 0$ such that 
$$	\PP \left( \langle Z_k , Z_{k+1} \rangle_o \ge \alpha k \right) \tends{ k \to \infty} 1 \ .$$

The proof of the well-definiteness of the escape rate is an almost direct application of Kingman's ergodic subadditive theorem. Indeed, it is a classical fact from dynamical systems that $T$ is ergodic, see for example \cite[Proposition 3.2]{coudi}. In order to be able to apply Kingman's ergodic subadditive theorem, one needs to show that $d(o,Z_n)$ is (up to an additive constant) a subadditive cocycle which follows from Lemma \ref{lemme.comparisondistance}. Using successively the triangular inequality and Lemma \ref{lemme.comparisondistance} we have for any $m,n \in \NN$
	\begin{eqnarray*}
		 d(o,Z_{n +m}) & \le & d(o,Z_{n}) + d(Z_n,Z_{n+m}) \\
		 & \le & d(o,Z_{n}) + d(o,Z_{m} \circ T^n) + R  \ . \
	\end{eqnarray*}

Consider the function $f_n := d(o, Z_n) +R$. By Lemma \ref{lemme.comparisondistance} the function  $f_n$ is subbaditive (ie.  $ f_{n+m} \leq f_m \circ T^n + f_n$). Since as we mentioned earlier, $T$ is ergodic, using Kingmann ergodic theorem, this implies that for almost all $\omega \in \Omega:$
$$\lim_{n\tv \infty} \frac{f_n(\omega)}{n}= \ell.$$
And we have $\frac{f_n(\omega)}{n} = \frac{d(o,Z_n) -R}
{n}.$ Passing to the limit proves that for almost all trajectories: 
$$\lim_{n\tv \infty} \frac{d(o,Z_n)}{n} = \ell \ .$$

\hfill $\blacksquare$

\subsection{Proof of Theorem \ref{th-main-asymptotic-behaviour-of-Zn}.}  This section is dedicated to the proof of \ref{th-main-asymptotic-behaviour-of-Zn}. The proof follows from a result of Mathieu and Sisto in \cite{artmathieusisto} that we explain below. \\

Let $\cF_n$ be a $\sigma$-field generated and $f_n$ be measurable with respect to $\cF_n$. We say in this case that $f_n$ is  a \textbf{defective adapted cocycle}, see \cite[Definition 3.1]{artmathieusisto}.
The \textbf{defect} of $F=(f_n)_{n\in \N}$ is the map $\Psi:=(\Psi_{n,m})_{(n,m)\in \N\times \N}$ defined by 
\begin{equation}
\Psi_{n,m}(\omega) = f_{n+m}(\omega) - f_n(\omega) - f_m(T_n \omega) \ . 
\end{equation}

A defective adapted cocycle $F=(f_n)$ is said to have a $p-th$ \textbf{finite moment} if, $\bE(|f_1|^p)$ is finite. \\ 

A defective adapted cocycle $F=(f_n)$ is said to satisfy a \textbf{second moment deviation inequality} if there is a constant $C$ such that for any $n,m \in \NN$, $ \EE(\Psi_{n,m}^2) \le C $. \\

In order to prove Theorem \ref{th-main-asymptotic-behaviour-of-Zn} we rely on the following theorem. We denote by $\bV(f)=\bE(f^2)-\bE(f)^2$ the variance of a random variable $f$.

\begin{theorem}\cite[Lemma 3.4, Theorem 4.1, Theorem 4.2]{artmathieusisto}
\label{theo_mathieusisto}
Let $F=(f_n)$ be a defective adapted cocycle. Assume that $F$ has a finite second moment and satisfies the second moment deviation inequality, then 
\begin{enumerate}
\item There exists $\tau\geq 0$ such that for all $n\in \N$, $\left| \frac{1}{n} \bE[f_n] -\ell\right| \leq \frac{\tau}{n}$.
\item There exists $\sigma\geq 0$ such that the  sequence $\frac{V(f_n)}{n}$ converges to $\sigma^2$. 
\item The law of $\frac{1}{\sqrt{n}} (f_n - \ell n ) $ weakly converges to the Gaussian law with zero mean and variance $\sigma^2$. 
\end{enumerate}
\end{theorem}

Theorem \ref{th-main-asymptotic-behaviour-of-Zn} follows from the above theorem applied to the following defective adapted cocyle. \\

Let $\cF_n$ be the $\sigma$-field generated by $X_n$. 
The function  $f_n=d(o,Z_n)$ is measurable with respect to $\cF_n$ and therefore  $d(o,Z_n)$ is  a defective adapted cocycle. Since $d(o,Z_1)$ has finite support it has in particular a finite exponential moment. In order to use the above theorem, one has then to verify that $F$ satisfies a second moment deviation inequality. \\

In this case, the defect of $F=(f_n)_{n\in \N}$ is given by 
\begin{equation}
\Psi_{n,m} =  d(o,Z_n) + d(o,Z_n \circ T_m) - d(o, Z_{m+n})  \ . 
\end{equation}

We will actually prove a stronger statement that the one required in order to apply Theorem \ref{theo_mathieusisto}: we will see that $F$ satisfies a exponential moment deviation inequalities.
 
\begin{proposition}
\label{prop_expdevinequalitiesdefec}
There is $\epsilon > 0$ such that for any $n,m \in \NN$ and any $R > 0$ we have
$$\bP( |\Psi_{n,m}| \geq  R ) = \epsilon^{-1} e^{-\epsilon R} \ . $$ 
\end{proposition} 

which will conclude the proof of Theorem \ref{th-main-asymptotic-behaviour-of-Zn}. The proof follows the same line as the one of \cite{artmathieusisto} but must be adapted to to our setting: in \cite{artmathieusisto} the authors study the pushforward of random walks on the isometry group of a hyperbolic space. \\

\textbf{Proof of Proposition \ref{prop_expdevinequalitiesdefec}} In  a first step, we reduce the proof of Proposition \ref{prop_expdevinequalitiesdefec} to the one of Proposition \ref{prop.morseprobabilistic} stated in the introduction. In fact, we show that $1/2 \cdot \Psi_{n,m}$ is up to a bounded error equal  to the Gromov product 
$$ \langle o, Z_{n+m} \rangle_{Z_n} :=  \frac{1}{2} \left( d(o,Z_n) + d(Z_n,Z_{n+m}) -d(o,Z_{n+m}) \right) \ .$$
This readily follows from Lemma \ref{lemme.comparisondistance} whose conclusion is indeed 
$$ |d(o,Z_n \circ T_m) -d(Z_n,Z_{n+m}) | \le R \ , $$ where $R>0$ is the diameter of $M_0$. Therefore, for any $n,m \in \NN$ 
\begin{equation}  
 \label{eq.Psi is at bounded distance of gromov product}
 	|\Psi_{n,m} - \langle Z_n, Z_{n+m}\rangle_{o} |\leq R 
\end{equation}
 
It remains then to show that Proposition \ref{prop.morseprobabilistic} holds. We recall here its statement for the reader's convenience. \\

\textit{	Let $(Z_n)_{n \in \NN}$ be a sequence of random variables valued in a geodesic Gromov hyperbolic space such that:}
		\begin{itemize}
			\item \textit{there is $R > 0$ such that for any $n \in \NN$, $d(Z_n, Z_{n+1}) \le R$}
			\item \textit{the sequence $(Z_n)_{n \in \NN}$ satisfies the linear progress with exponential tail property.}
		\end{itemize}
	\textit{Then it satisfies an exponential deviation inequality: there is $\epsilon > 0$ such that for any $ n,m \in \NN$ and any $t > 0$ we have}
$$\bP(\left<Z_{n+m},Z_0 \right>_{Z_n}  \geq  t ) = \epsilon^{-1} e^{-\epsilon t} \ . $$ 

\textbf{Proof of Proposition} \ref{prop.morseprobabilistic}. Now we follow the strategy of Mathieu-Sisto \cite[Section 11]{artmathieusisto}. The main geometric ingredient is the following lemma, established in any Gromov $\delta$-hyperbolic space: 
\begin{lemma}\cite[Lemma 11.4]{artmathieusisto}\label{lem-geometric_lemma}
For all $\epsilon>0$ there exists $C>0$ with the following property. Let $(w_i)_{i\in [0,n]}$ be a sequence of points and denote by $\g$ the geodesic between $w_0$ and $w_n$. For any $T\geq C$, for any $k\in [0, n]$ one of the following holds:
\begin{enumerate}
\item[H1.] There exist $k_1<k\leq k_2$ with $|k_2-k_1|\leq T$ so that $d(w_{k_i}, w_{k_i+1})\geq (d(w_k,\g) -C)/T$  for $i\in \{ 1, 2\}$.
\item[H2.] There exist $k_1<k\leq k_2$ with $|k_2-k_1|\geq T$ so that $d([w_{k_1}, w_{k_1+1}], [w_{k_2}, w_{k_2+1}]) \leq \epsilon (k_2-k_1)$.
\item[H3.] There exist $k_1<k\leq k_2$ with $|k_2-k_1|\geq T$ so that $\sum_{i \in [k_1, k_2)} d(w_i, w_{i+1}) \geq e^{(k_2-k_1)/C}/C$.
\end{enumerate}
\end{lemma}

We want to control the probability that $\langle o, Z_{n}\rangle_{Z_k} \geq t$ occurs. Note that it is enough to control $\langle o, Z_{n}\rangle_{Z_k} \geq t$ for $t \ge t_0$ for some fixed $t_0$. \\

We are going to apply this lemma to $w_i = Z_i(\omega)$. In particular, since $d(Z_i,Z_{i+1}) \leq R$, we will be able to eliminate case H1 and H3 by choosing $R$ large enough. We will then see how satisfying H2 implies the conclusion of Proposition \ref{prop.morseprobabilistic}. We start by adjusting the constants $\epsilon, t_0$ and $C$ in order for H1 not to occur. \\

We start by $\epsilon$. Recall that from Proposition  \ref{prop-linearprogress} that there exists  $\alpha >0$  such that for all $k_1,k_2\geq 1$:
\begin{equation}\label{spectral -equation}
\bP(d(Z_{k_1},Z_{k_2}) \leq \alpha  |k_1-k_2| ) \leq \frac{1}{\alpha} e^{-\alpha|k_1-k_2|}.
\end{equation} 
We set $\epsilon = \frac{\alpha}{2}$. Let $C>0$ be the corresponding constant given by Lemma \ref{lem-geometric_lemma}. \\

Before adjusting $t_0$, we set 
$$ T = T(t) := (t -\delta -C)/(R+1) \ , $$
where $\delta$ is the hyperbolicity constant. We now adjust $t_0$ large enough such that for any $t \ge t_0$ we have
 \begin{equation}
 	\label{eqdef.toetT}
   \left\{ \begin{array}{l}
		  T(t) \ge C \\
		  R C \ T(t) < e^{T(t)/C} \\
		  T(t) \ge t / [2(R+1)] \ ,  
	\end{array}  \right. 
\end{equation}
which is possible since $T$ is linear in $t$.  \\

Let us now see why H1 does not occur under then event $\langle o, Z_{n}\rangle_{Z_k} \ge t \ge t_0$. By $\delta$-hyperbolicity, we have
$$ \langle o, Z_{n}\rangle_{Z_k} \geq t \Longrightarrow d(\omega_k, \gamma) \ge t - \delta \ .$$

Note that in the H1 case, one must have
	\begin{eqnarray*}
		 R & \ge & d(Z_{k_i}, Z_{k_i+1})   \ge (d(Z_k,\g) -C)/T \ge (t - \delta -C)/T \\ 
		 & \ge & (t - \delta -C)/T \ , 
	 \end{eqnarray*}
which contradicts the first of the inequalities of \eqref{eqdef.toetT}. \\

Let us now see why H3 does not occur either. Let $|k_1 -k_2| > T(t)$ as in H3. Note that the triangular inequality gives
 $$\sum_{i \in [k_1, k_2)} d(Z_i, Z_{i+1})\leq |k_2-k_1|R \ , $$ 
 
 which is smaller than  $e^{(k_2-k_1)/C}/C$ by the second inequality of \eqref{eqdef.toetT}. \\

Hence H2 must occur: for $t \ge t_0$ large enough,  if $\langle o, Z_{n}\rangle_{Z_k} \geq t$,  there exists $k_1< k \leq k_2$ with $|k_2-k_1|\geq T(t)$ so that 
$$d([Z_{k_1}, Z_{k_1+1}], [Z_{k_2}, Z_{k_2+1}]) \leq \epsilon (k_2-k_1).$$

Using the triangle inequality and the fact that $\alpha = 2 \epsilon$ we get 
$$ d(Z_{k_1}, Z_{k_2}) \leq\frac{\alpha}{2}(k_2-k_1)+ 2R \leq \alpha (k_2-k_1) \ . $$

 Therefore, for any $t \ge t_0$ there are $k_2 - k_1 \ge T(t)$ such that 
	$$ \PP(\langle o, Z_{n}\rangle_{Z_k} \geq t) \le   \PP(d(Z_{k_1}, Z_{k_2}) \leq \alpha (k_2-k_1)). $$
We now use the linear progress with exponential tail property to get 
	\begin{eqnarray*}
		\PP(\langle o, Z_{n}\rangle_{Z_k} \geq t) & \le  & \alpha^{-1} e^{- \alpha (k_2 -k_1)} 
	 \end{eqnarray*}
 
We sum the above inequality on the set $ \{ (k_1,k_2) \in [0,n]^2)\, | \, |k_1-k_2|\geq T(t)\} \  $
by pairing the pairs $(k_2,k_1)$ according to the value of the difference $k_2 -k_1$. Indeed, we have
 $$\sharp \{ (k_2,k_1) \ |  k_1 \le k \le k_2 \ \text{ and } \ k_2 - k_1 = p \} \le p^2 \ .$$
And then
$$ \PP(\langle o, Z_{n}\rangle_{Z_k} \geq t)  \le  \alpha^{-1} \somme{ p \ge T(t)} p^2 \  e^{- \alpha p } \ .$$
Therefore
	$$ \PP(\langle o, Z_{n}\rangle_{Z_k} \geq t)  \le   C_2 \ e^{- \alpha T(t) } \ ,$$

for some constant $C_2$. We conclude by using the third inequality of \eqref{eqdef.toetT}. \hfill $\blacksquare$ \\

\appendix

\section{Dirichlet and standard random walks.}
\label{sec.dirichletstandard} 

This section aims at investigating when a Dirichlet random walk is actually a standard random walk. We start by showing that this is the case for flat torus and their covers. We will then show that it is the only essentially non compact example.

\subsection{The flat torus case.} We explain in this subsection how the study of the Dirichlet random walk in the case of the universal cover of a flat torus reduces to a simple random walk on $\R^k$. Intermediate coverings could be treated in the same way. \\

Let $T^d$ be a flat torus of dimension $d$. Up to rescaling the metric, one can assume that $\mu_g(T^d) =1$. We denote by $dx$ the Lebesgue measure. We fix our starting point $o := 0_{\RR^d}$ (and then $p_0 = 0_{T^d}$). Let $\G \subset \R^d$ be such that $T^d$ is isometric to $\R^d/\G$. The group $\G$ is isomorphic to $\Z^d$ and acts by translations on $\R^d$. Let $\tau_x$ be the translation of vector $x\in \R^k$. \\

The reason that the Dirichlet random walk behaves like a usual random walk is that Dirichlet domains 'commute' with translations. Indeed, we have that for all $x,y \in \R^d$, 
	$$\tau_y D_x = D_{\tau_y x}= D_{y+x}. $$
	
Indeed, let $z\in \tau_y D_x$. We have for all $\g \in \G$,
$$d(x,  \tau^{-1}_y z) \leq d(x, \g \tau^{-1}_y z) \ . $$ 

Since $\G$ is a group of translations, we see that $\g$ is acting as $\tau_v$ for some $v\in \R^d$. In particular 
$$\g \tau^{-1}_y = \tau_v \tau^{-1}_y  = \tau^{-1}_y \tau_v = \tau^{-1}_y \g \ . $$

Therefore for all $\g\in \G$ we have
$d(\tau_y x, z) \leq d(\tau_y x, \g z) \ .$ Let $d\mu (x) := 1_{D_o} (x) dx$. The above discussion implies in particular that we have 
	\begin{equation}
		\label{eq.toruscase}
		 d [ ({\tau_y})^*(\mu)] (x) = 1_{D_o}(x-y)dx = 1_{D_{(\tau_y o)}}(x)dx  \ .
	\end{equation}
Let $Y_n$ be a sequence of I.I.D. random variables with distribution $\mu$: $Z_n$  follows the law of $S_n:= \sum_{k=1}^n Y_k$. All the classical results for random walks on $\R^d$ therefore apply in the case of the Dirichlet random walk. For example, Polya's Theorem implies that the sequence $Z_n$ on $\R^d$ is transcient if and only if $d \geq 3$. In the next subsection, we will see that it is essentially the only occurrence of such a phenomenon.

\subsection{When is a Dirichlet random walk standard.}

This subsection is devoted to classify Dirichlet random walks that are standard random walks. Recall that a standard random walk is defined as the pushforward of a random walk of the isometry group of $M$. More precisely, we say that a stochastic process $(Z_n)_{n \in \NN}$ of $M$ is a \textbf{standard random walk} if there is a measure $\mu$ on $G := \Isom(X)$ and a point $o \in M$ such that for any $n \in \NN$, $Z_n$ follows the law of $ \nu_1 \cdot ... \cdot \nu_n \cdot o $ where $(\nu_i)_{i \in \NN}$ are independent and identically distributed according to $\mu$. \\

Note first that the Dirichlet random walk can visit any point on $M$. Therefore, the group $G_{\mu}$ generated by the support of $\mu$ must act transitively on $M$, so $M$ has to be a homogeneous space. This is already a far stronger assumption that the one we assumed in all of the theorems stated in the introduction, as for example, it implies that $M$ has constant scalar curvature. One can then reduce this question  to homogeneous space and turn it into a question of Lie group theory. \\

The example of the previous subsection shows that a Dirichlet random walk can be a standard random walk. Note also that one can construct compact examples. For example by considering the case of a round $n$-projective space, the quotient of a round $n$-sphere by the antipodal map. Both these examples are actually of the same nature: the deck group commutes with a subgroup of isometries acting transitively on the covering. In the second case, the antipodal map is in the center of $O(n)$, and in the second case the full group of isometries is itself Abelian. Let us now come to the following simple characterisation of such examples. 

\begin{proposition}
\label{prop.dirichletegalstandard}
	Let $M$ be a complete connected Riemannian manifold such that the Dirichlet random walk is a standard random walk. Then $M$ is a metric product of $\mathbb{R}^d$ with a compact homogeneous space.
\end{proposition}

\textbf{Proof.} We start with introducing some notations. Let $\Gamma$ the cocompact lattice of $G := \Isom(M)$ such that $M_0 = \quotient{M}{\Gamma}$. Let $\overline{G_{\mu}}^0$ be the connected component containing $\Id$ of the closure in $G$ of the group $G_{\mu}$ generated by $\supp{\mu}$. By Cartan's closed subgroup theorem, $\overline{G_{\mu}}^0$ is a Lie group. Since $G_{\mu} \cdot o = M$, we have as well that $\overline{G_{\mu}}^0 \cdot o = M$ (this set is open and closed and $M$ is connected). \\

The key ingredient of the proof is the following lemma.

\begin{lemma}
\label{lemma.dirichletegalstandard}
Under the assumption of Proposition \ref{prop.dirichletegalstandard} and with the notations introduced above. For any $\nu \in \overline{G_{\mu}}^0$ and any $ \gamma \in \Gamma$ we have $\gamma \nu = \nu \gamma$.
\end{lemma}

Before proving the above lemma, let us see how it implies the conclusion of Proposition \ref{prop.dirichletegalstandard}. \\

\textbf{Proof of (Lemma  \ref{lemma.dirichletegalstandard} $\Rightarrow$ \textbf{Proposition} \ref{prop.dirichletegalstandard})}. Let $G' < G$ be the subgroup generated by both $\overline{G_{\mu}}^0$ and $\Gamma$. Note first that the mapping 
$$ \fonctionbis{\overline{G_{\mu}}^0 \times \Gamma}{G'}{(g, \gamma)}{ g\cdot \gamma } $$
is a morphism since $\Gamma$ and $\overline{G_{\mu}}^0$ commute. Note also that $H :=  \overline{G_{\mu}}^0 \cap \Gamma$ is Abelian and in the center of both $\Gamma$ and $\overline{G_{\mu}}^0$. In particular we have
$$ G' = \quotient{(\overline{G_{\mu}}^0 \times \Gamma)}{H} \ , $$
where $H \curvearrowright \overline{G_{\mu}}^0 \times \Gamma$ by $h \cdot (g, \gamma) = ( h \cdot g, h^{-1} \cdot \gamma)$. Let us see that  this implies that $H$ is cocompact in both $\Gamma$ and $\overline{G_{\mu}}^0$. Indeed, since both $\overline{G_{\mu}}^0$ and $\Gamma$ are cocompact in $G$ (the first one because $\overline{G_{\mu}}^0 \cdot o = M$ and second one by assumption) they must be cocompact in $G'$ as well. In particular
	$$ \quotient{G'}{\overline{G_{\mu}}^0} = \quotient{\Gamma}{H} \hspace{0.2cm} \text{ and } \hspace{0.2cm} \quotient{G'}{\Gamma} = \quotient{\overline{G_{\mu}}^0}{H} $$
are compact. Therefore, up to a finite index, $\Gamma$ must be isomorphic to $\ZZ^d$ as $H$ is Abelian. This implies that $\overline{G_{\mu}}^0$ is the product of an Abelian group and a compact Lie group. We conclude by recalling once again that $\overline{G_{\mu}}^0$ acts transitively on $M$.  \hfill $\blacksquare$ \\

Let us now prove that Lemma \ref{lemma.dirichletegalstandard} holds. \\

\textbf{Proof of Lemma \ref{lemma.dirichletegalstandard}.} We shall actually show that for all $\gamma \in \Gamma$ and all $h \in \overline{G_{\mu}}^0$ we have 
	\begin{equation}
		\label{eq.gmuandgammacommute}
		 \gamma \cdot h \cdot o = h \cdot \gamma \cdot o \ . 
	\end{equation}

This is equivalent to the statement above. Indeed, since $\overline{G_{\mu}}^0 \cdot o = M$ we have $h_x \in \overline{G_{\mu}}^0$ such that $h_x \cdot o = x$ for any $x \in M$. Therefore, for any $h \in \overline{G_{\mu}}^0$, any $\gamma \in \Gamma$ and any $x \in M$
\begin{eqnarray*}
	\gamma \cdot h \cdot x &=& \gamma \cdot h \cdot h_x  \cdot o \\
	&=& h \cdot (h_x \cdot \gamma \cdot o) \\
	&=& h \cdot (\gamma \cdot h_x \cdot o) \\
	&=& h \cdot \gamma \cdot x \ ,
\end{eqnarray*}
using repeatedly \eqref{eq.gmuandgammacommute}. Let us now focus on proving that \eqref{eq.gmuandgammacommute} holds. \\

Recall that assuming that the Dirichlet random walk is a standard random walk corresponds to assuming that $Z_n$ follows the law of $\nu_1 \cdot ... \cdot \nu_n \cdot o$ where the $\nu_i$s are I.I.D random variables following the law of some measure $\mu$. \\

If $x \in M$ and $y \in M_0$, recall that we denoted by $D_{x}(y)$ the unique lift of $y$ in the $D_x$. Since we assumed that $Z_n = \nu_1 \cdot ... \cdot \nu_n \cdot o$ in law, one has, by definition of the Dirichlet random walk, the following equality in law 
$$ D_{\nu_1 \cdot ... \cdot \nu_n \cdot o}(X_{n+1}) = \nu_1 \cdot ... \cdot \nu_{n+1} \cdot o \ . $$
Thefore	
$$  (\nu_1 \cdot ... \cdot \nu_{n})^{-1} D_{\nu_1 \cdot ... \cdot \nu_n \cdot o}(X_{n+1}) = \nu_{n+1} \cdot o \ . $$
Since $\nu_{n+1}$ is independent of $\nu_1, ..., \nu_{n}$ and since it follows the law of $\mu$, we have for any $n$, $\mu^{\otimes n}$-almost surely
$$ (\nu_1 \cdot ... \cdot \nu_{n})^{-1} D(\nu_1 \cdot ... \cdot \nu_n  \cdot o) = D(o) \ . $$
Note that any open set of the closure $\overline{G_{\mu}}^0$ of $G_{\mu}$ must be charged by $\mu^{*N}$ for some $N \in \NN$. In particular, it implies by continuity that for all $\nu \in \overline{G_{\mu}}^0$ one has 
	\begin{equation}
		\label{eq.dirichletstandard2}
			 \nu^{-1} D(\nu \cdot o) = D(o) \ . 	
	\end{equation}

Let us now see that the above equality implies the conclusion of Proposition \ref{prop.dirichletegalstandard}. Let $S$ be the generating set of $\Gamma$ given by the elements defining the Dirichlet domain $D(o)$: the elements $s \in \Gamma$ such that there is $x \in D(o) $ such that 
	$$ d(x, s \cdot o) = \mini{\gamma \in \Gamma} \  d(x, \gamma \cdot o) \ .$$
In other words, more geometric, these are the elements that pair faces of the Dirichlet domain $D(o)$.
	
	 We show that for any $s \in S$ and any $\nu$ in  small neighbourhood $U$ of the identity of $\overline{G_{\mu}}^0$ we have 	
\begin{equation}	
	\label{eq.dirichletstandard}
	 \nu^{-1} \cdot s \cdot \nu \cdot o = s \cdot o \ , 
\end{equation}
which implies Proposition \ref{prop.dirichletegalstandard} since $S$ (resp. $U$) generates $\Gamma$ (resp. $\overline{G_{\mu}}^0$). Let us then show that \eqref{eq.dirichletstandard} holds. \\

We start by a general observation. For any compact set $K$ of $M$ that contains $o$, for any $\gamma \in G$ (in particular in $\Gamma$) and any $\epsilon > 0$ there is a neighbourhood $U$ of the identity in $G$ (in particular in $\overline{G_{\mu}}^0$) such that for any $\nu \in U$ and any $x \in K$
	$$ d(x, \gamma \cdot o) - \epsilon \le d(x, \nu^{—1} \gamma \nu \cdot o) \le d(x, \gamma \cdot o) + \epsilon \ .$$
The proof of the upper bound consist in using successively the triangular inequlity:
\begin{eqnarray*}
 d( \nu^{—1} \gamma \nu \cdot o,x) & \le &  d( \nu^{—1} \gamma \nu \cdot o , \nu^{-1} \cdot x) + d( \nu^{-1} \cdot x, x) \\
 & \le &  d( \gamma \nu \cdot o, \gamma \cdot o) + d(\gamma \cdot o,x) + d( \nu^{-1} \cdot x, x) \\
 & \le &  d( \nu \cdot o, o) + d(\gamma \cdot o,x) + d( \nu^{-1} \cdot x, x) \ , 
\end{eqnarray*}
which is smaller than $\epsilon$ whenever $\supr{x \in K} \ d( \nu \cdot x, x) \le \epsilon/2$. The lower bound follows from the same argument. \\

The above remark applied with $K = D(o)$ together with the fact that the $\Gamma \cdot o$ is discrete shows that there is a neighbourhood $U$ of the identity in $\overline{G_{\mu}}^0$ such that for any $\nu \in U$ we have 
$$ \nu^{-1} D(\nu \cdot o) = \{ x \in M \ , \ d(x, o) \le \mini{s \in S} \ d( x, \nu^{-1} s \nu \cdot o) \}.$$

We argue by contradiction and assume that there is $s \in S$ such that $s \cdot o \neq \nu^{-1} s \nu \cdot o$ under the assumption that $\nu^{-1} D(\nu \cdot o) = D(o)$. Since we assume that $\nu \in U$, we know a priori that $  \nu^{-1} s \nu \cdot o $ is close to $s \cdot o$ (in particular it cannot be any other point of $S \cdot o$). \\

\begin{figure}[h!]
\begin{center}
	\def\svgwidth{0.4 \columnwidth}
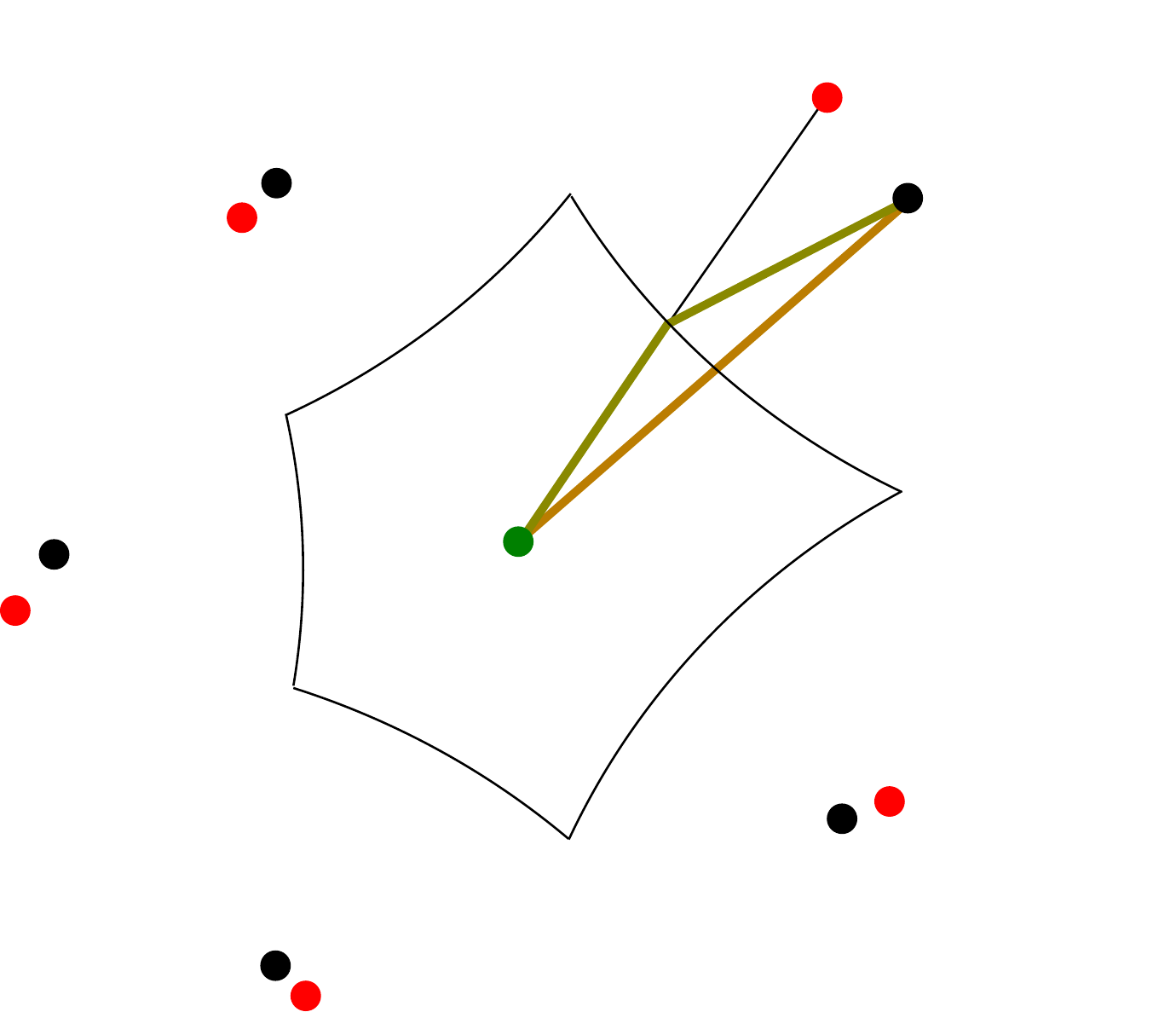
	\end{center}
\caption{The red dots correspond to the set $\{ \nu^{-1} s \nu \cdot o \ , \ s \in S \}$ and the black ones to $S \cdot o$.}
	\label{fig.polygon}
\end{figure}

We consider a geodesics from $o$ to $\nu^{-1} s \nu \cdot o$. Such a geodesic must intersect the median plane  $ \{ z \in M \ , \ d(z,o) = d(z, s \cdot o) \}$ at a different point that any of the geodesics relating $o$ and $s \cdot o$. In particular, by construction and using the absurd assumption, one could then be able to construct a broken minimising geodesic from $o$ to $s \cdot o$ (in green in Figure \ref{fig.polygon}). This contradicts that geodesics are smooth. \hfill $\blacksquare$

\bibliographystyle{alpha}

\end{document}

%% file: preambule_spectral_gap.tex
\usepackage{geometry}
\usepackage[T1]{fontenc}    
\usepackage[utf8]{inputenc} 
\usepackage{amssymb,amsmath,amsthm}

\newcommand{\R}{\mathbb{R}}

\newcommand{\N}{\mathbb{N}}
\newcommand{\Z}{\mathbb{Z}}

\newcommand{\bE}{\mathbb{E}}

\newcommand{\bP}{\mathbb{P}}

\newcommand{\bV}{\mathbb{V}}

\newcommand{\cA}{\mathcal{A}}

\newcommand{\cF}{\mathcal{F}}

\newcommand{\cO}{\mathcal{O}}

\newcommand{\G}{\Gamma}
\newcommand{\g}{\gamma}

\newcommand{\tv}{\rightarrow}

\DeclareMathOperator{\Isom}{Isom}

\DeclareMathOperator{\Id}{Id}

\DeclareMathOperator{\Vol}{Vol}

\newtheorem{theorem}{Theorem}[section]
\newtheorem{definition}[theorem]{Definition}
\newtheorem{lemma}[theorem]{Lemma}
\newtheorem{proposition}[theorem]{Proposition}

\usepackage[english]{babel}
\usepackage{ucs}
\usepackage{tikz}
\usepackage{amsmath}
\usepackage{amsthm}
\usepackage{amsfonts}
\usepackage{amssymb}
\usepackage{hyperref}
\usepackage{dsfont}
\usepackage{mathrsfs}
\usepackage{eqnarray}
\write18{}
\usepackage{xcolor}
\usepackage{graphicx}
\usepackage{pdfpages} 
\usepackage[nothm]{thmbox}
\usepackage{esint}

\hypersetup{
    colorlinks=true,
    linkcolor=red,
    filecolor=magenta,      
    urlcolor=cyan,
}

\numberwithin{equation}{section}

\newcommand{\EE}{\mathbb{E}}

\newcommand{\PP}{\mathbb{P}}
\newcommand{\HH}{\mathbb{H}}
\newcommand{\ZZ}{\mathbb{Z}}
\newcommand{\RR}{\mathbb{R}}

\newcommand{\NN}{\mathbb{N}}
\newcommand{\inte}[1]{\underset{#1}{\int}}

\newcommand{\supp}[1]{\mathrm{supp}(#1)}

\newcommand{\limiinf}[1]{\underset{#1}{\liminf}}
\newcommand{\mini}[1]{\underset{#1}{\min}}

\newcommand{\infi}[1]{\underset{#1}{\inf}}
\newcommand{\supr}[1]{\underset{#1}{\sup}}
\newcommand{\somme}[1]{\underset{#1}{\sum}}

\newcommand{\tends}[1]{\underset{#1}{\longrightarrow}}

\newcommand{\quotient}[2]{{\raisebox{.2em}{$#1$}\left/\raisebox{-.2em}{$#2$}\right.}}

\newcommand{\fonction}[5]{\begin{array}{r r c l}
					#1 \hspace{2mm} : & #2 & \to & #3 \\
					& #4 & \mapsto & #5 \\
			   \end{array}}
\newcommand{\fonctionbis}[4]{\begin{array}{r c l}
					   #1 & \to & #2 \\
					 #3 & \mapsto & #4 \\
			   \end{array}}

\DeclareMathOperator{\vol}{vol}

\theoremstyle{definition}

\newtheorem{remark}[equation]{Remark}

\theoremstyle{theorem}

\newtheorem{theoreme}[equation]{Theorem}

%% file: dirrandomwalk.pdf_tex
\begingroup%
  \makeatletter%
  \providecommand\color[2][]{%
    \errmessage{(Inkscape) Color is used for the text in Inkscape, but the package 'color.sty' is not loaded}%
    \renewcommand\color[2][]{}%
  }%
  \providecommand\transparent[1]{%
    \errmessage{(Inkscape) Transparency is used (non-zero) for the text in Inkscape, but the package 'transparent.sty' is not loaded}%
    \renewcommand\transparent[1]{}%
  }%
  \providecommand\rotatebox[2]{#2}%
  \newcommand*\fsize{\dimexpr\f@size pt\relax}%
  \newcommand*\lineheight[1]{\fontsize{\fsize}{#1\fsize}\selectfont}%
  \ifx\svgwidth\undefined%
    \setlength{\unitlength}{747.60673702bp}%
    \ifx\svgscale\undefined%
      \relax%
    \else%
      \setlength{\unitlength}{\unitlength * \real{\svgscale}}%
    \fi%
  \else%
    \setlength{\unitlength}{\svgwidth}%
  \fi%
  \global\let\svgwidth\undefined%
  \global\let\svgscale\undefined%
  \makeatother%
  \begin{picture}(1,0.43178894)%
    \lineheight{1}%
    \setlength\tabcolsep{0pt}%
    \put(0,0){\includegraphics[width=\unitlength,page=1]{dirrandomwalk.pdf}}%
    \put(0.20233637,0.02861177){\color[rgb]{0.64313725,0,0}\makebox(0,0)[lt]{\lineheight{1.25}\smash{\begin{tabular}[t]{l}$X_o = p_o$\end{tabular}}}}%
    \put(0.70654759,0.00707333){\color[rgb]{0.65490196,0,0}\makebox(0,0)[lt]{\lineheight{1.25}\smash{\begin{tabular}[t]{l}$o$\end{tabular}}}}%
    \put(0.15547335,0.25317172){\color[rgb]{0,0.52941176,0}\makebox(0,0)[lt]{\lineheight{1.25}\smash{\begin{tabular}[t]{l}$X_1$\end{tabular}}}}%
    \put(0.34848508,0.12615931){\color[rgb]{0,0.52941176,0}\makebox(0,0)[lt]{\lineheight{1.25}\smash{\begin{tabular}[t]{l}$X_2$\end{tabular}}}}%
    \put(0.03531356,0.10654322){\color[rgb]{0,0.52941176,0}\makebox(0,0)[lt]{\lineheight{1.25}\smash{\begin{tabular}[t]{l}$X_3$\end{tabular}}}}%
    \put(0.83545318,0.08299768){\color[rgb]{0,0.52941176,0}\makebox(0,0)[lt]{\lineheight{1.25}\smash{\begin{tabular}[t]{l}$Z_1$\end{tabular}}}}%
    \put(0.73284066,0.19147378){\color[rgb]{0,0.52941176,0}\makebox(0,0)[lt]{\lineheight{1.25}\smash{\begin{tabular}[t]{l}$Z_2$\end{tabular}}}}%
    \put(0.87311112,0.40090925){\color[rgb]{0,0.52941176,0}\makebox(0,0)[lt]{\lineheight{1.25}\smash{\begin{tabular}[t]{l}$Z_3$\end{tabular}}}}%
    \put(0,0){\includegraphics[width=\unitlength,page=2]{dirrandomwalk.pdf}}%
  \end{picture}%
\endgroup%

%% file: dirichlet_domain.pdf_tex
\begingroup%
  \makeatletter%
  \providecommand\color[2][]{%
    \errmessage{(Inkscape) Color is used for the text in Inkscape, but the package 'color.sty' is not loaded}%
    \renewcommand\color[2][]{}%
  }%
  \providecommand\transparent[1]{%
    \errmessage{(Inkscape) Transparency is used (non-zero) for the text in Inkscape, but the package 'transparent.sty' is not loaded}%
    \renewcommand\transparent[1]{}%
  }%
  \providecommand\rotatebox[2]{#2}%
  \newcommand*\fsize{\dimexpr\f@size pt\relax}%
  \newcommand*\lineheight[1]{\fontsize{\fsize}{#1\fsize}\selectfont}%
  \ifx\svgwidth\undefined%
    \setlength{\unitlength}{841.88976378bp}%
    \ifx\svgscale\undefined%
      \relax%
    \else%
      \setlength{\unitlength}{\unitlength * \real{\svgscale}}%
    \fi%
  \else%
    \setlength{\unitlength}{\svgwidth}%
  \fi%
  \global\let\svgwidth\undefined%
  \global\let\svgscale\undefined%
  \makeatother%
  \begin{picture}(1,0.70707071)%
    \lineheight{1}%
    \setlength\tabcolsep{0pt}%
    \put(0,0){\includegraphics[width=\unitlength,page=1]{dirichlet_domain.pdf}}%
    \put(0.31218892,0.57094181){\color[rgb]{0,0,0}\makebox(0,0)[lt]{\lineheight{1.25}\smash{\begin{tabular}[t]{l}$\Gamma \cdot x$\end{tabular}}}}%
    \put(0.49913929,0.3277248){\color[rgb]{1,0,0}\makebox(0,0)[lt]{\lineheight{1.25}\smash{\begin{tabular}[t]{l}$x$\end{tabular}}}}%
    \put(0.47009842,0.26419802){\color[rgb]{0,0.54901961,0}\makebox(0,0)[lt]{\lineheight{1.25}\smash{\begin{tabular}[t]{l}$D_x$\end{tabular}}}}%
  \end{picture}%
\endgroup%

%% file: selffat.pdf_tex
\begingroup%
  \makeatletter%
  \providecommand\color[2][]{%
    \errmessage{(Inkscape) Color is used for the text in Inkscape, but the package 'color.sty' is not loaded}%
    \renewcommand\color[2][]{}%
  }%
  \providecommand\transparent[1]{%
    \errmessage{(Inkscape) Transparency is used (non-zero) for the text in Inkscape, but the package 'transparent.sty' is not loaded}%
    \renewcommand\transparent[1]{}%
  }%
  \providecommand\rotatebox[2]{#2}%
  \newcommand*\fsize{\dimexpr\f@size pt\relax}%
  \newcommand*\lineheight[1]{\fontsize{\fsize}{#1\fsize}\selectfont}%
  \ifx\svgwidth\undefined%
    \setlength{\unitlength}{1028.46287746bp}%
    \ifx\svgscale\undefined%
      \relax%
    \else%
      \setlength{\unitlength}{\unitlength * \real{\svgscale}}%
    \fi%
  \else%
    \setlength{\unitlength}{\svgwidth}%
  \fi%
  \global\let\svgwidth\undefined%
  \global\let\svgscale\undefined%
  \makeatother%
  \begin{picture}(1,0.44664132)%
    \lineheight{1}%
    \setlength\tabcolsep{0pt}%
    \put(0,0){\includegraphics[width=\unitlength,page=1]{selffat.pdf}}%
  \end{picture}%
\endgroup%

%% file: almostsubadditive.pdf_tex
\begingroup%
  \makeatletter%
  \providecommand\color[2][]{%
    \errmessage{(Inkscape) Color is used for the text in Inkscape, but the package 'color.sty' is not loaded}%
    \renewcommand\color[2][]{}%
  }%
  \providecommand\transparent[1]{%
    \errmessage{(Inkscape) Transparency is used (non-zero) for the text in Inkscape, but the package 'transparent.sty' is not loaded}%
    \renewcommand\transparent[1]{}%
  }%
  \providecommand\rotatebox[2]{#2}%
  \newcommand*\fsize{\dimexpr\f@size pt\relax}%
  \newcommand*\lineheight[1]{\fontsize{\fsize}{#1\fsize}\selectfont}%
  \ifx\svgwidth\undefined%
    \setlength{\unitlength}{651.15739981bp}%
    \ifx\svgscale\undefined%
      \relax%
    \else%
      \setlength{\unitlength}{\unitlength * \real{\svgscale}}%
    \fi%
  \else%
    \setlength{\unitlength}{\svgwidth}%
  \fi%
  \global\let\svgwidth\undefined%
  \global\let\svgscale\undefined%
  \makeatother%
  \begin{picture}(1,0.79896029)%
    \lineheight{1}%
    \setlength\tabcolsep{0pt}%
    \put(0,0){\includegraphics[width=\unitlength,page=1]{almostsubadditive.pdf}}%
    \put(0.27254148,0.72151906){\color[rgb]{1,0,0}\makebox(0,0)[lt]{\lineheight{1.25}\smash{\begin{tabular}[t]{l}$\gamma \cdot o$\end{tabular}}}}%
    \put(0.13408592,0.14892324){\color[rgb]{1,0,0}\makebox(0,0)[lt]{\lineheight{1.25}\smash{\begin{tabular}[t]{l}$o$\end{tabular}}}}%
    \put(0.92045561,0.4704217){\color[rgb]{0.60784314,0.50196078,0}\makebox(0,0)[lt]{\lineheight{1.25}\smash{\begin{tabular}[t]{l}$Z_{n +m}$\end{tabular}}}}%
    \put(0.82423004,0.00812103){\color[rgb]{1,0.58431373,0}\makebox(0,0)[lt]{\lineheight{1.25}\smash{\begin{tabular}[t]{l}$Z_m \circ T^n$\end{tabular}}}}%
    \put(0.36640959,0.46572834){\color[rgb]{0.00392157,0.50196078,0}\makebox(0,0)[lt]{\lineheight{1.25}\smash{\begin{tabular}[t]{l}$Z_n$\end{tabular}}}}%
    \put(0.41334367,0.64877123){\color[rgb]{0.60784314,0.50196078,0}\makebox(0,0)[lt]{\lineheight{1.25}\smash{\begin{tabular}[t]{l}$Z_{n+1}$\end{tabular}}}}%
    \put(0.26080793,0.11137607){\color[rgb]{1,0.58431373,0}\makebox(0,0)[lt]{\lineheight{1.25}\smash{\begin{tabular}[t]{l}$Z_1 \circ T^n$\end{tabular}}}}%
    \put(0,0){\includegraphics[width=\unitlength,page=2]{almostsubadditive.pdf}}%
    \put(0.14112603,0.54316956){\color[rgb]{0,0,0}\makebox(0,0)[lt]{\lineheight{1.25}\smash{\begin{tabular}[t]{l}$\gamma^{-1}$\end{tabular}}}}%
  \end{picture}%
\endgroup%

%% file: polygons.pdf_tex
\begingroup%
  \makeatletter%
  \providecommand\color[2][]{%
    \errmessage{(Inkscape) Color is used for the text in Inkscape, but the package 'color.sty' is not loaded}%
    \renewcommand\color[2][]{}%
  }%
  \providecommand\transparent[1]{%
    \errmessage{(Inkscape) Transparency is used (non-zero) for the text in Inkscape, but the package 'transparent.sty' is not loaded}%
    \renewcommand\transparent[1]{}%
  }%
  \providecommand\rotatebox[2]{#2}%
  \newcommand*\fsize{\dimexpr\f@size pt\relax}%
  \newcommand*\lineheight[1]{\fontsize{\fsize}{#1\fsize}\selectfont}%
  \ifx\svgwidth\undefined%
    \setlength{\unitlength}{397.37402581bp}%
    \ifx\svgscale\undefined%
      \relax%
    \else%
      \setlength{\unitlength}{\unitlength * \real{\svgscale}}%
    \fi%
  \else%
    \setlength{\unitlength}{\svgwidth}%
  \fi%
  \global\let\svgwidth\undefined%
  \global\let\svgscale\undefined%
  \makeatother%
  \begin{picture}(1,0.85958353)%
    \lineheight{1}%
    \setlength\tabcolsep{0pt}%
    \put(0,0){\includegraphics[width=\unitlength,page=1]{polygons.pdf}}%
    \put(0.81261726,0.67712786){\color[rgb]{0,0,0}\makebox(0,0)[lt]{\lineheight{1.25}\smash{\begin{tabular}[t]{l}$s \cdot o$\end{tabular}}}}%
    \put(0.72201285,0.8014874){\color[rgb]{1,0,0}\makebox(0,0)[lt]{\lineheight{1.25}\smash{\begin{tabular}[t]{l}$\nu^{-1} s \nu \cdot o$\end{tabular}}}}%
    \put(0.39873293,0.30819453){\color[rgb]{0,0.45098039,0}\makebox(0,0)[lt]{\lineheight{1.25}\smash{\begin{tabular}[t]{l}$o$\end{tabular}}}}%
  \end{picture}%
\endgroup%